\documentclass[12pt]{article}

\usepackage{amsmath,amsfonts,amssymb}
\usepackage[all]{xy}
\usepackage{graphicx}
\usepackage{epsf}
\usepackage{cancel}
\usepackage{multirow}
\usepackage{rotating}  
\usepackage{arydshln}
\usepackage{tikz}
\usepackage{eucal}
\usepackage{enumitem}
\usepackage{nicefrac}

\usetikzlibrary{decorations.markings}
\usetikzlibrary{arrows}
\usepackage[colorlinks=true,urlcolor=black]{hyperref}

\hypersetup{			
backref = true,			
pagebackref = true,		
hyperindex = true, 		
colorlinks = true, 		
breaklinks = true, 		
urlcolor = blue, 		
linkcolor = blue, 		
bookmarks = true,		
bookmarksopen = true,  
citecolor=red,
}

\usepackage{hyperref}
\usepackage{breakurl}
\hypersetup{breaklinks=true}

\advance \textheight by \topskip
\newtheorem{theorem}{Theorem}[section]

\newtheorem{corollary}[theorem]{Corollary}

\newtheorem{definition}[theorem]{Definition}
\newtheorem{example}[theorem]{Example}

\newtheorem{lemma}[theorem]{Lemma}

\newtheorem{proposition}[theorem]{Proposition}
\newtheorem{remark}[theorem]{Remark}

\newtheorem{DP}[theorem]{Definition/Proposition}

\newcommand{\qeed}{\hfill\textrm{QED}\break\null}
\newenvironment{demo}{\noindent\textit{Proof.}~}{\qeed}

\def\F{{\cal F}}

\def\FF{\mathbb F}

\def\g{\mathfrak{g}}

\newcommand{\beq}{\begin{equation}}
\newcommand{\eeq}{\end{equation}}
\newcommand{\beqa}{\begin{eqnarray}}
\newcommand{\eeqa}{\end{eqnarray}}
\newcommand{\A}{{\mathrm {Aut}}}

\newcommand{\noi}{\noindent}
\newcommand{\nn}{\nonumber}

\newcommand{\e}{\epsilon}
\def\>{\rangle}
\def\<{\langle}

\begin{document}

\title{Commutation relations of $\g_2$  and the incidence geometry of the Fano plane}
%\begin{flushright}
%{\small USACH-FM-01-02}\\[-0.4cm]
%{\small PM-01-07}\\[1cm]
%\end{flushright}
{\bf }

\author{
{\sf M. Rausch de Traubenberg}\thanks{e-mail:
michel.rausch@iphc.cnrs.fr}$\,\,$${}^{a}$ and
{\sf M. J. Slupinski}\thanks{e-mail:
marcus.slupinski@math.unistra.fr}$\,\,$${}^{b}$
\\
{\small ${}^{a}${\it IPHC-DRS, UdS, CNRS, IN2P3 
}}\\
{\small {\it  23  rue du Loess, Strasbourg, 67037 Cedex, France 
}}\\
{\small ${}^{b}${\it
Institut de Recherches en Math\'emathique Avanc\'ee,
UdS and CNRS}}\\
{\small {\it 7 rue R. Descartes, 67084 Strasbourg Cedex, France.}}  \\ 
%{\small ${}^{c}${\it Laboratoire de Math\'ematiques et Applications,
%Universit\'e de Haute Alsace,}} \\
%{\small {\it Facult\'e des Sciences et Techniques,
%4 rue des Fr\`eres Lumi\`ere, F-68093 Mulhouse, France.}}
}

\maketitle
\date
\vskip-1.5cm

\vspace{2truecm}

\begin{abstract}
\noindent
We continue our
study  and classification
of  structures on the Fano plane $\F$ and its dual $\F^\ast$
involved  in the construction of octonions and the Lie algebra $\g_2(\FF)$ over a %n arbitrary
field
$\FF$.
These
are a ``composition factor''  $\epsilon : \F\times \F \to \{-1,1\}$,hhhh
inducing an octonion multiplication, and a function $\delta^\ast :  \A(\F)\times \F^\ast \to \{-1,1\}$
such that  $g\in \A(\F)$  can be lifted to  an automorphism of the octonions {\it iff} $\delta^\ast(g,\cdot)$
is the Radon transform of a function on $\F$.
We
lift the action of Aut$(\F)$ on $\F$ to the action  of a non-trivial eight-fold covering  Aut$(\hat \F_\epsilon)$  on a two-fold covering  $\hat \F_\epsilon $ of $\F$
contained in the octonions. This 
extends tautologically  to an action  on the octonions by automorphism. 
Finally, we  associate to incident  point-line pairs  a generating set   of $\g_2(\FF)$
and express  brackets  in terms of the incidence geometry of $\F$  and
$\epsilon$.
\end{abstract}

\vspace{2cm}

\newpage
We continue our study and classification of structures on the Fano plane ${\cal F}$ and its dual  ${\cal F}^\ast$   involved in the construction of octonions and the Lie algebra $\mathfrak g_2 (\mathbb F)$ over a field $\mathbb F$. These are a "composition factor" :  ${\cal F}\times {\cal F} \to\{-1, 1\}$, inducing an octonion multiplication, and a function $\delta^\ast : Aut({\cal F}) \times {\cal F}^\ast \to \{-1, 1\}$ such that $g \in  Aut({\cal F})$ can be lifted to an automorphism of the octonions iff $\delta^\ast(g, \cdot)$ is the Radon transform of a function on ${\cal F}$. We lift the action of $Aut({\cal F})$  on ${\cal F}$ to the action of a non-trivial eight-fold covering   $Aut({\cal F})$ on a twofold covering $\hat {\cal F}$ of ${\cal F}$  contained in the octonions. This extends tautologically to an action on the octonions by automorphism. Finally, we associate to incident point-line pairs a generating set of $\mathfrak g_2 (\mathbb F)$ and express brackets in terms of the incidence geometry of ${\cal F}$ and .

\newpage
%\section{Introduction}
In 1900 F. Engel \cite{Eng} realised the exceptional complex Lie group $G_2(\mathbb C)$ as the linear isotropy group  of a generic three-form in seven dimensions, and in 1907 his student W. Reichel \cite{rei}
realised the real compact exceptional Lie group $G^c_2(\mathbb R)$ as the linear isotropy group  of a real three-form in seven dimensions.
\'E. Cartan  \cite{Car} showed that the automorphism group of the octonions is the exceptional Lie group $G^c_2(\mathbb R)$  and, to the best of the authors' knowledge, 
 H. Freudenthal \cite{freu} was the first to point out  a close relationship between the Fano plane $\F$ and
the octonions.
For a more detailed historical perspective see \cite{agri}.
In this article we give a systematic
study of  additional structures  on $\F$, over and above the projective  structure, which are relevant   to the  construction of   the octonions
over any field  $\FF$  of characteristic not two
 (see also \cite{mm}).
In particular we give  explicit generators  of $\g_2(\FF)$,   the Lie algebra   of $G^c_2(\FF)$,
 together with their brackets in terms of  ``augmented'' incidence relations of $\F$.
 We also give simple formul\ae \ for
the  $\g_2(\FF)-$invariant three and four  forms  on the space of imaginary octonions
in terms of these relations.
\\

Our starting point is Aut$(\F)$, the group of automorphisms of $\F$. This is a simple group of order $168$ isomorphic to
$GL(3,\mathbb Z_2)$ or $PSL(2,\mathbb Z_7)$ \cite{atlas}.
Elements of order seven, which we consider as ``orientations'' of $\F$, fall into 
 two conjugacy classes which we  consider as defining two   ``pre-orientations'' of $\F$.
Two pre-orientations are distinguished by the description
of lines in $\F$ in terms of any compatible orientation.

To an orientation of $\F$ one can associate a composition factor on $\F$.
A  composition factor
is the  structure allowing a generalisation  of Freudenthal's construction of a composition algebra product on an eight-dimensional vector space $\mathbb O_\F$ canonically associated to $\F$.
Our first result completes the classification of composition factors given in \cite{mm}. 
  We  show that for the action of $\A(\F)$ on  composition factors  there are two orbits, each containing eight elements.
Given a composition factor $\epsilon$ on $\F$ it is then natural to ask whether the action of Aut$(\F)$ on $\F$ extends to an action by automorphisms on 
$\mathbb O_\F$ equipped with the  product corresponding to $\epsilon$.
This is not the case but  it  turns out that this action  lifts to an action of an eight-fold non-split
extension of Aut$(\F)$  \cite{wil}.
In order  to understand this phenomenon  we first 
associate to  $g\in$Aut$(\F)$ a function $\delta^\ast(g,\cdot)$ on $\F^\ast$, the space of lines in $\F$.
We then show that solving the lifting problem for  $g$  is equivalent to
finding a function on $\F$ whose Radon transform is $\delta^\ast(g,\cdot)$ and in this way solve it.

In the last part of the paper we associate to  each incident pair $(P,D) \in \F\times \F^\ast$ an element  $X_{P,D}$ of $\g_2(\FF)$.
This gives twenty-one elements  which span the fourteen-dimensional vector space $\g_2(\FF)$. We express the brackets of
the $X_{P,D}$  in terms of the incidence relations of $\F$  and the composition factor $\epsilon$.  Finally, this allows us to identify various geometric subalgebras of $\g_2(\FF)$ associated
to  points and  lines of $\F$. In particular, to each point we  associate a Cartan subalgebra and an
$\mathfrak{su}(3)$(or $\mathfrak{sl}(3)$)  subalgebra  containing it. Dually, we associate 
 to each line  an $\mathfrak{so}(3)$ and an
$\mathfrak{so}(3) \times \mathfrak{so}(3)$  subalgebra containing it as an ideal.\\

We now give a more detailed description of the contents of this paper.
In Section \ref{sec:F}
we recall the definition and basic properties of the Fano plane $\F$
and its dual $\F^\ast$. 
 We introduce the notion of
an orientation of $\F$ and 
 show that such an orientation  induces an orientation of each line and of each point.
In Section \ref{sec:cv-19}, having recalled relevant definitions and results from \cite{mm}, we first 
show that the action of $\A(\F)$ on composition factors has two orbits. We then consider the Radon transform from functions on $\F$ to functions
on $\F^\ast$, and identify its kernel and image. The exponential version of the Radon transform turns out to be very
useful in the context of this article.  
Fixing a composition factor $\epsilon$
the next question we address is whether the obvious linear action of $g \in$ Aut$(\F)$ on
the associated composition algebra $(\mathbb O_\F,\epsilon)$
is by automorphism. A necessary and sufficient condition for this is that for any line $D$ of $\F$, $g$ induces an isomorphism
 of the subalgebras associated to $D$ and $g \cdot D$.
The function  $\delta^\ast:\A(\F) \times \F^\ast \to  \{-1,1\}$ which detects this property is introduced and for fixed $g$ classified
(Theorem \ref{theo:delta*}): there are eight possibilities, each of
which is realised by twenty-one  elements of Aut$(\F)$. 
In particular the  linear action of $g\in\A(\F)$ on $(\mathbb O_\F,\epsilon)$ is by automorphism {\it iff}  $\delta^\ast(g,P)=1, \forall P
\in \F$.
To conclude Section \ref{sec:cv-19} we introduce a ``double covering''  $\hat \F_\epsilon$  of the pair  $(\F,\epsilon)$
(see also  p. 205  in \cite{KT}). The group of automorphisms of $(\mathbb O_\F,\epsilon)$ is canonically isomorphic to the
group of automorphisms of $\hat \F_\epsilon$ by restriction, and
we show that $\pi: \A(\hat \F_\epsilon) \to \A(\F)$ is a non-split  extension
 by $\mathbb Z_2^3$ (see also \cite{wil, KT}).
The main point here is to show that  finding a lift of   $g\in \A(\F)$ that acts by automorphism on  $\hat \F_\epsilon$
is equivalent to finding a function on $\F$ whose Radon transform is $\delta^\ast(g,\cdot)$.

In Section \ref{sec:g2} we fix a composition factor and our
 starting point is the observation that, notwithstanding  non-associativity,
left multiplication by  purely imaginary octonions gives
a representation  of the Clifford algebra of Im$(\mathbb O_\F)$  acting on $\mathbb O_\F$ (see for example \cite{wdm}).
This means one can realise the Lie algebra $\mathfrak g_2$ as the annihilator in $\mathfrak{so}(\text{Im}(\mathbb O_\F))$  of $1 \in \mathbb O_\F$ \cite{har}. To each incident pair $(P,D) \in \F \times \F^\ast$ we associate an element $X_{P,D}$ of $\g_2(\FF)$ realised in this way.
In order to describe the action of $\A(\hat \F_\epsilon)$ on the $X_{P,D}$ we introduce a function $\delta : \A(\hat \F_\epsilon) \times \F \to \{-1,1\}$.
As $\hat g$ varies in Aut$(\hat \F_\epsilon)$ this gives rise to sixty-four functions on $\F$,
each of which is realised by twenty-one elements of $\A(\hat \F_\epsilon)$.
It is a surprising fact  that
the exponential of the Radon transform of $\delta(\hat g, \cdot)$ is exactly $\delta^\ast(\pi(g), \cdot)$ for 
 $\hat g \in \A(\hat\F_\epsilon)$.
To get explicit formul\ae \ for  brackets of the $X_{P,D}$ we find simple normal forms for  the orbits of the action of $\A(\hat \F_\epsilon)$ on
pairs of  incident pairs  $\big((P,D),(P',D')\big)  \in \big(\F \times \F^\ast\big)^2$.
These formul\ae \ are  Fano geometric in the sense that they only involve the incidence geometry of $\F$ and  the composition factor (Theorem \ref{theo:g2CR}).
The point  here is that it easy to guess these brackets up to  a sign but rather more subtle to determine the signs precisely.  
Furthermore,  using the formul\ae \   we associate   to each point of $\F$ a Cartan subalgebra and an $\mathfrak{su}(3)$ (or $\mathfrak{sl}(3))$  subalgebra containing it.
Dually, to each line in $\F$ we  associate an $\mathfrak{so}(3)$ subalgebra and  an $\mathfrak{so}(3)\times\mathfrak{so}(3)$  subalgebra containing it as an ideal. In this way
any rank two subalgebra of
$\g_2(\FF)$ can be realised (up to conjugation) as a subalgebra associated either to a line or to  a point.\\

%\begin{enumerate}
%\item[]
\begin{minipage}{15.85cm}
{\it Throughout this paper: $\FF$  denotes a field of characteristic not two,  $\mathbb Z_n$ denotes the ring $\mathbb Z/n \mathbb Z$
and  $S_2$ denotes the group $\{-1,1\}$. The elements of $\mathbb Z_2$ will be denoted $0,1$ (not $[0], [1]$) when there is no ambiguity.}
%\end{enumerate}
\end{minipage}
\section{Properties of the Fano plane}\label{sec:F}
In this section we recall without proof  the basic properties of the Fano plane we need in the rest of the paper.
\subsection{The Fano plane and its dual}
A projective plane is a set of points together with a collection of subsets called lines such that
two distinct lines intersect in a unique point and two distinct points are contained in a unique line.
Given a projective plane ${\cal P}$ the dual plane ${\cal P}^*$ is defined as 
 the set whose points are the lines of ${\cal P}$,
and whose lines are the subsets of ${\cal P}^*$ consisting
of  concurrent lines. It is easily checked that  ${\cal P}^*$ is also a projective plane.

If $V$ is a three-dimensional vector space over a field $\mathbb F$, ${ P}(V)$, the
set of lines passing through the origin in $V$ together with the collection of subsets consisting  of lines contained in a fixed plane
passing through the origin of $V$, defines a projective plane.

\begin{figure}[!ht]
  %\begin{center} \begin{tikzpicture}[scale=1.2]
  \begin{center}

\vskip -1truecm
    
  \begin{tikzpicture}[scale=.8] 
\tikzstyle{point}=[circle,draw]

\tikzstyle{ligne}=[thick]
\tikzstyle{pointille}=[thick,dotted]

\tikzset{->-/.style={decoration={
  markings,
  mark=at position .5 with {\arrow{>}}},postaction={decorate}}}

\tikzset{middlearrow/.style={
        decoration={markings,
            mark= at position 0.5 with {\arrow{#1}} ,
        },
        postaction={decorate}
    }}

\coordinate  (3) at ( -2,-1.15);%   {$\bullet$}; %[point] {$e_3$};%(-a/2,-asqrt(3)/6)

\coordinate  (2) at  ( 0,-1.15) ;% {$\bullet$};%[point] {$e_2$};%(0,-asqrt(3)/6)
\coordinate  (5) at   ( 2,-1.15);% {$\bullet$};% [point] {$e_5$};%(a,-asqrt(3)/6)

\coordinate  (6) at   (0,2.31);% {$\bullet$};% [point] {$e_6$};%(0,a(sqrt(3)/2-sqrt(3)/6)
\coordinate  (1) at  (1, .58);% {$\bullet$};% [point] {$e_1$};%(a/4,a(sqrt(3)/4-sqrt(3)/6)

\coordinate  (4) at  (-1., .58) ;% {$\bullet$}; %[point] {$e_4$};
\coordinate  (7) at  (0,0) ;%{$\bullet$};% [point] {$e_7$};

% -> latex oriente la droite
\draw  (6) --(1);
\draw (1) --(5);

\draw  (5) --(2);
\draw  (2) --(3);

\draw  (3) --(4);
\draw  (4) --(6);

\draw  (6) --(7);
\draw  (7) --(2);

\draw (3) --(7);
\draw  (7) --(1);

\draw (5) --(7);
\draw  (7) --(4);

\draw (1) to [bend left=65] (2);
\draw  (2) to [bend left=65] (4);
%\draw [middlearrow={triangle 45}] (4) to [bend left=65] (1);
\draw  (4) to [bend left=65] (1);
%\draw (1)  node[scale=0.6] {$\bullet$} ;
\draw (1)  node {$\bullet$} ;

\draw (2)  node{$\bullet$} ;

\draw (3)  node{$\bullet$} ;
\draw (4)  node{$\bullet$} ;

\draw (5)  node{$\bullet$} ;

\draw (6)  node{$\bullet$} ;

\draw (7)  node{$\bullet$} ;

  \end{tikzpicture} \hskip 3.truecm 
%\begin{figure}%[!ht]
  %\begin{center}
    %\begin{tikzpicture}[thick,scale=.8]
    \begin{tikzpicture}[thick,scale=.6]
  \coordinate (P0) at (0,0);
    \coordinate (P1) at (4.5,0) ; %(3,0);
    \coordinate (P4) at (6.6,2.1); %(5.1,2.1); 
    \coordinate (P2) at (2.1,2.1);

     \coordinate (P3) at (0,3);
    \coordinate (P7) at  (4.5,3) ;% (3,3);
    \coordinate (P6) at (6.6,5.1); %  (5.1,5.1); 
    \coordinate (P5) at (2.1,5.1);

 \draw[solid][line width=1pt] (P0) -- (P1);
 \draw[solid][line width=1pt] (P0) -- (P2);
  \draw[solid][line width=1pt] (P0) -- (P3);
  \draw[solid][line width=1pt] (P0) -- (P4);
   \draw[solid][line width=1pt] (P0) -- (P5);
    \draw[solid][line width=1pt] (P0) -- (P6);
 
      \draw[solid][line width=1pt] (P0) -- (P7);

\draw[dashed][line width=1pt] (P0) -- (P1)-- (P4) -- (P2)--(P0);
    
    \draw[dashed][line width=1pt] (P3) -- (P7)-- (P6) -- (P5)--(P3);
    \draw[dashed][line width=1pt] (P0) --(P3);
    \draw[dashed][line width=1pt] (P1) --(P7);
    \draw[dashed][line width=1pt] (P4) --(P6);
    \draw[dashed][line width=1pt] (P2) --(P5);

\draw (P0)  node {$\bullet$} ;
\draw (P1)  node {$\bullet$} ;

\draw (P2)  node{$\bullet$} ;

\draw (P3)  node{$\bullet$} ;
\draw (P4)  node{$\bullet$} ;

\draw (P5)  node{$\bullet$} ;

\draw (P6)  node{$\bullet$} ;

\draw (P7)  node{$\bullet$} ;

%\draw (P1) node[left]{$1$} ;
%    \draw (P2) node[ right]{$2$} ;
%    \draw (P3) node[above left]{$3$} ;
%    \draw (P4) node[above right]{$4$} ;
%    \draw (P5) node[ left]{$5$} ;
%    \draw (P6) node[right]{$6$} ;
%    \draw (P7) node[right=2pt]{$7$} ;
   \draw (P0) node[below left=2pt]{$O$} ;

   \end{tikzpicture}
  
\caption{The Fano plane ${\cal F}$ and the Fano cube $V_{\cal F}$.}
\label{fig:FP}
 \end{center}
\end{figure}
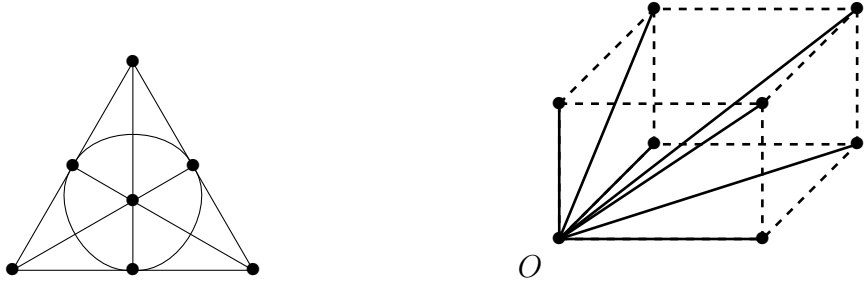

The simplest example of a projective plane %is
 %the Fano plane ${\cal F}$ which
 consists  of seven points and seven lines arranged  as below,
%%%%%% fano-plane
%\begin{minipage}{7cm}
%\input{fano-plane}
%\end{minipage}
%%%%%%%%%%%%
%%%% fano-cube
%\begin{minipage}{7cm}
%\input{fano-cube}
%\end{minipage}
%%%%%%%%%
each line contains  exactly three points and each point is contained in exactly  three lines (see Figure \ref{fig:FP}).\\
\vskip .7truecm

Any projective plane containing seven points and seven lines can be  obtained by projectivising a canonically associated vector space:
\begin{definition}
Let $\F$ be a projective plane  consisting of seven points and seven lines. Let $ V_\F= \F \cup \{0\}$ be the set obtained by
formally adding  a point $0$ to $\F$. Define the symmetric map $+:  V_\F \times  V_\F \to  V_\F$ by
\beqa
P+Q=\left\{
\begin{array}{cll}
R & \text{~if~} & P \ne Q \text{~and~} P,Q,R \text{~are~aligned}\ ;\\
0 & \text{~if~} & P = Q\ ; \\
P & \text{~if~} & Q = 0 \ .
\end{array}\right.\nn
\eeqa
With respect to the obvious scalar multiplication by $\mathbb Z_2$ this defines the structure of a three-dimensional
$\mathbb Z_2-$vector space on $ V_\F$  with zero element $0$.  The natural bijection $P (V_\F) \cong \F$ is an isomorphism of projective planes and
if $P, Q, R$ are three non-zero points in $V_\F$, then $P, Q, R$  are coplanar {\it iff} $P+Q+R=0$.

The set of lines of $\F$ will be denoted $\F^\ast$. This is a projective plane whose `lines' are triples of concurrent lines so in particular
three lines  $D_1,D_2,D_3$ are concurrent iff  $D_1+D_2+D_3=0$. 
\end{definition}
Any two  three-dimensional vector spaces over $\mathbb Z_2$ are isomorphic so 
any two projective planes consisting of seven points and seven lines are isomorphic in
the sense that there is a bijection between them which sends lines to lines.
From now on we will refer to
any projective plane consisting of seven lines and seven points as a Fano plane 
and  any three-dimensional vector space over $\mathbb Z_2$ as a Fano cube.
\subsection{Automorphisms of the Fano plane}
%%%%%%%%%%%%%%%%%%%%%%%%%%%%%%%%%%%%%%%%%%%%%%%%%%%%%%%%%%%%%%%%%%%%%%%%%%%%%%%%%%%%%%

\begin{definition}
Let $\F$ be  a Fano plane. We set 
\beqa
\mathrm{Aut}(\F) =
\big\{ f : \F \to \F \ \ \mathrm{s.t. ~} f \mathrm{~ is ~a~bijection
~and~} f \mathrm{~sends~lines~to~lines}\big\} \ .\nn
\eeqa
 This group is also called the group of collineations of $\F$.
\end{definition}

Any bijection  $f$:  $\F \to \F$ extends to a unique bijection $\hat f:$
$V_\F\setminus\{0\}\to V_\F\setminus\{0\}$, and clearly this defines  a natural group isomorphism
\beqa
\text{Aut}(\F) \cong GL\big(V_\F\big) \ .\nn
\eeqa
It is well-known  \cite{atlas} that this general linear group is simple of order $168$,  and generated by two elements $a$ and $b$ satisfying the relations 
\beqa
\label{eq:ab}
a^2=b^3=(ab)^3=(aba^{-1}b^{-1})^4=1 \ . 
\eeqa
%One can also show that it is isomorphic to
%\beqa
%\text{Aut}(\FF) \cong GL_2(\mathbb Z_7) \ . \nn
%\eeqa

We now summarise the main properties of $\A(\F)$.

\begin{proposition}\label{cor:4}{\color{white} toto} 
\begin{enumerate}
\item Let $f \in \A(\F)$ such that $f \ne 1$. Then $f$ is of order
$2,3,4$ or $7$.
\item Let $f\in \A(\F)$ be of order two.
Then there exist a line $L\in \F^\ast$ and a point $P\in\F$ such that $f(R)=R, \forall R\in L$ and $R+f(R) +P=0$.
 $\forall R\not \in L$,
Conversely, given $L \in \F^\ast$ and $P\in L$ there exists a unique $f\in\A(\F)$ of order two such that $f(R)=R,
\forall R\in L$ and $R+f(R)+P=0, \forall R\not \in L$. Two elements of order two are conjugate.
\item Let $f\in \A(\F)$ be of order three. Then there
exists a unique  triangle  stable by $f$. Conversely every triangle
is obtained in this way from exactly two order three elements of
$\A(\F)$. Two elements of order three are conjugate.
\item  Let $f\in \A(\F)$ be of order four.  Then there exists $P\in\F$ and $L\in\F^\ast$ containing $P$ such that
 $f(P)=P$,  $f(L)=L$ and
$f$ is of order two on $L$, and  for all $R \not \in L$, $R+f(R)+P\ne 0$. Conversely, given  $P\in \F$ and $L\in\F^\ast$ containing
$P$, there exists exactly two elements of order four such that $f(P)=P$ and $f(L)=L$. Two elements of order four are conjugate.
\item   Let $f\in GL\big(V_\F\big)$ be of order seven. Then the  minimal polynomial of $f$ is either $x^3+x^2+1$ or $x^3+x+1$.
Two elements of order seven are conjugate iff they have the same minimal polnomial.
\end{enumerate}
\end{proposition}
\begin{corollary}
In $GL\big(V_\F\big)$ there are respectively $21, 56$, $42$ and $48$  elements respectively  of order respectively two, three, four and seven.
\end{corollary}

Recall that if $n \in \mathbb N$ the Legendre symbol  $\Big( \frac n 7 \Big)$ can be  defined by:
$
\Big( \frac n 7 \Big)=  n^3 \mod 7\ .  
$

\begin{corollary}\label{cor:legendre}
Let $f \in GL\big(V_\F\big)$ be of order seven and let $n,m$ be  positive integers. Then
$f^n$ is conjugate to $f^m$  {\it iff}   $\Big( \frac n 7 \Big)=\Big( \frac m 7 \Big)$.
\end{corollary}

\begin{example}\label{ex:ab}
Recall that according to  \eqref{eq:ab} Aut$(\F)$ can be generated by $a,b$ satisfying $a^2=b^3=(ab)^7=(aba^{-1}b^{-1})^4=1$.
Denoting the points of $\F$ temporarily by the numbers $1,\cdots,7$, an
explicit example (see Figure \ref{fig:OF}) is given by
\beqa
a= \begin{pmatrix} 1&2&3&4&5&6&7\\
                   1&2&7&4&6&5&3
                   \end{pmatrix} \ , \ \
 b= \begin{pmatrix} 1&2&3&4&5&6&7\\
                    2&7&4&6&5&3&1
                   \end{pmatrix}   \ .\nn               
\eeqa
It the follows that
\beqa
ab=\begin{pmatrix} 1&2&3&4&5&6&7\\
                   2&3&4&5&6&7&1
                   \end{pmatrix}=\tau \ , \ \
aba^{-1}b^{-1}=   \begin{pmatrix} 1&2&3&4&5&6&7\\
                                4&2&6&1&7&5&3
                   \end{pmatrix}          \ . \nn        
\eeqa 

\end{example}
\subsection{{Orientation, orientation type and lines}}\label{sec:orientation}
\begin{definition}
{\color{white} t}

\begin{enumerate}
\item Let $\F$ be a Fano plane. An orientation of $\F$ is an element   $\tau \in \text{Aut}(\F)$ of order seven.
\item Let $(\F,\tau)$ and $(\F',\tau')$ be oriented Fano planes. We say $(\F,\tau)$ and $(\F',\tau')$
are isomorphic
{\it iff} there exists an isomorphism of Fano planes  $f :\F\rightarrow\F'$ s.t. $f \circ \tau \circ f^{-1}=\tau'$.
\end{enumerate}
\end{definition}
{ It is well-known that the lines of a Fano plane can be described in one of  two ways with respect to an orientation.
\begin{proposition}\label{prop:droites}
Let $(\F,\tau)$ be an oriented Fano plane. Then
we have one of the following:
\begin{enumerate}
\item[(i)] for all $P$ in $\F$, the triple $D_P=\Big\{P,\tau(P),\tau^3(P)\Big\}$ is a line and each line can be
written uniquely in this way;
\item[(ii)] for all $P$ in $\F$, the triple $D_P=\Big\{P,\tau^2(P),\tau^3(P)\Big\}$ is a line and each line can be
written uniquely in this way.
\end{enumerate}
\end{proposition}

%\begin{demo}
%Let $ V_\F$ be the Fano cube corresponding to $\F$ (see Proposition \ref{prop:cube})
%and $\hat \tau \in GL( V_\F)$ be the lift of $\tau$ to $V_\F$.
%To prove the Proposition it is sufficient to prove that for each $P\ne 0$ in $V_\F$, the planes of $V_\F$
%can be written uniquely as 
%\beqa
%\{0,P,\hat \tau(P), \hat \tau^3(P)\}\ , \nn
%\eeqa
%or as
%\beqa
%\{0,P,\hat \tau^2(P), \hat \tau^3(P)\}\ .\nn
%\eeqa
%By Proposition \ref{prop:auto} the minimal polynomial of $\hat \tau$ is either $x^3+x+1$ or
%$x^3+x^2+1$. In the first case this means that for any $P\ne 0$ in $ V_\F$ we have
%\beqa
%P+ \hat\tau(P) + \hat \tau^3(P) = 0 \ .\nn
%\eeqa
%However, $(P,\hat \tau(P), \hat \tau^3(P))$ are three distinct points ($\hat \tau^7=1$ and $7$ is prime) and hence by Corollary
%\ref{cor:droite} they are coplanar. It is then straightforward to see that if $Q\ne P$  the planes
%$(0,P,\hat \tau(P), \hat \tau^3(P))$ and $(0,Q,\hat \tau(Q), \hat \tau^3(Q))$ are distinct. This proves part (i) of the proposition
%and part (ii) follows similarly. 
%\end{demo}

\begin{definition}\label{rem:type}
An oriented Fano plan $(\F,\tau)$ is of type $(0,1,3)$ if the lines are as in (i) above and of
type $(0,2,3)$ if the lines are as in (ii) above.
\end{definition}
Since there are only two conjugacy classes of elements of order seven in $\text{Aut}(\F)$ it follows that:
\begin{corollary}
Two oriented Fano planes  are isomorphic iff they are of the same  type.
\end{corollary}
If $(\F,\tau)$ is an oriented Fano plane, the induced map $\tau^\ast:\F^\ast\rightarrow \F^\ast$  is an orientation  of the dual Fano plane  $\F^\ast$. 
\begin{proposition}\label{cor:droites}
Let $(\F,\tau)$ be an oriented Fano plane and let $P \in \F$. Then the three lines containing $P$ are
 $D_P,D_{\tau^{-1}(P)}$ and $D_{\tau^{-3}(P)}$ in case (i) above, and 
 $D_P,D_{\tau^{-2}(P)}$ and $D_{\tau^{-3}(P)}$ in case (ii) above.
In particular,  $(\F,\tau)$ and $(\F^\ast,\tau^{\ast}{^{-1}})$ are isomorphic whereas  $(\F,\tau)$ and $(\F^\ast,\tau^{\ast})$ are not  isomorphic.
\end{proposition}

  An orientation $\tau$ of a Fano plane $\F$ induces an orientation ({\it i.e.}, an  order three bijection) on every line  $D$ in $\F$  as follows (see Figure  \ref{fig:OF}). If
  $\tau$ is of type $(0,1,3)$ (resp. $(0,2,3)$)  there is a unique $P\in \F$ such that $D=\Big\{P,\tau(P),\tau^3(P)\Big\}$ (resp. $D=\Big\{P,\tau^2(P),\tau^3(P)\Big\}$). The 
  orientation induced by $\tau$ on $D$ is then the cyclic permutation
 \beqa\label{eq:orP}
(P,\tau(P),\tau^3(P))\qquad(\text{resp. }(P,\tau^2(P),\tau^3(P))) \ . 
\eeqa
 %%%%%%%%%%%%%%%%%%%%
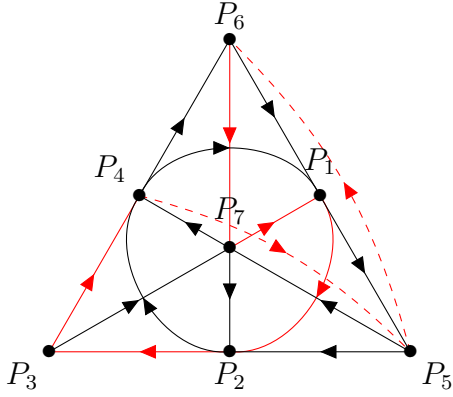
\begin{figure}[!ht]
\begin{center} \begin{tikzpicture}[scale=1.2]
\tikzstyle{point}=[circle,draw]

\tikzstyle{ligne}=[thick]
\tikzstyle{pointille}=[thick,dotted]

\tikzset{->-/.style={decoration={
  markings,
  mark=at position .5 with {\arrow{>}}},postaction={decorate}}}

\tikzset{middlearrow/.style={
        decoration={markings,
            mark= at position 0.5 with {\arrow{#1}} ,
        },
        postaction={decorate}
    }}

\coordinate  (3) at ( -2,-1.15);%   {$\bullet$}; %[point] {$e_3$};%(-a/2,-asqrt(3)/6)

\coordinate  (2) at ( 0,-1.15) ;% {$\bullet$};%[point] {$e_2$};%(0,-asqrt(3)/6)
\coordinate  (5) at ( 2,-1.15);% {$\bullet$};% [point] {$e_5$};%(a,-asqrt(3)/6)

\coordinate  (6) at (0,2.31);% {$\bullet$};% [point] {$e_6$};%(0,a(sqrt(3)/2-sqrt(3)/6)
\coordinate  (1) at (1, .58);% {$\bullet$};% [point] {$e_1$};%(a/4,a(sqrt(3)/4-sqrt(3)/6)

\coordinate  (4) at (-1., .58) ;% {$\bullet$}; %[point] {$e_4$};
\coordinate  (7) at (0,0) ;%{$\bullet$};% [point] {$e_7$};

\draw [color=red, middlearrow={triangle 45}] (1) to [bend left=65] (2);
\draw [color=red,  middlearrow={triangle 45}] (2) --(3);
\draw [color=red, middlearrow={triangle 45}] (3) --(4);
\draw [color=red,  middlearrow={triangle 45}] (6) --(7);
\draw [color=red,  middlearrow={triangle 45}] (7) --(1);
\draw[color=red,dashed,middlearrow={triangle 45}] (4) to [bend left=15]  (5);

\draw[color=red,dashed,middlearrow={triangle 45}] (5) to [bend right=15]  (6);
% -> latex oriente la droite
\draw [middlearrow={triangle 45}] (6) --(1);
\draw [ middlearrow={ triangle 45}] (1) --(5);

\draw[middlearrow={triangle 45}]  (5) --(2);

\draw[middlearrow={triangle 45}]  (3) --(7);

\draw  [middlearrow={triangle 45}]  (5) --(7);
\draw [middlearrow={triangle 45}] (7) --(4);
\draw[middlearrow={triangle 45}]  (4) --(6);
\draw [middlearrow={triangle 45}]  (7) --(2);

\draw [middlearrow={triangle 45}] (2) to [bend left=65] (4);
%\draw [middlearrow={triangle 45}] (4) to [bend left=65] (1);
\draw [middlearrow={triangle 45}]  (4) to [bend left=65] (1);
\draw (1)  node[scale=0.6] {$\bullet$} ;
\draw (1)  node {$\bullet$} ;
\draw (1) node[above=4pt  ]{$P_1$} ;

\draw (2)  node{$\bullet$} ;
\draw (2) node[below ]{$P_2$} ;
\draw (3) node[below left=0pt]{$P_3$} ;
\draw (3)  node{$\bullet$} ;
\draw (4)  node{$\bullet$} ;
\draw (4) node[ above left =1pt]{$P_4$} ;
\draw (5)  node{$\bullet$} ;
\draw (5) node[ below right ]{$P_5$} ;
\draw (6)  node{$\bullet$} ;
\draw (6) node[ above]{$P_6$} ;
\draw (7)  node{$\bullet$} ;
\draw (7) node[ above=5pt]{$P_7$} ;

%\draw[dashed,middlearrow={triangle 45}] (4) to [bend left=15]  (5);

%\draw[dashed,middlearrow={triangle 45}] (5) to [bend right=15]  (6);

%\draw[dashed,middlearrow={triangle 45}] (6) to [bend right=15]  (3);
%\draw[dashed,middlearrow={triangle 45}] (3) to [bend right=15]  (5);
%\draw[dashed,middlearrow={triangle 45}] (2) to [bend right=15]  (6);
%\draw[dashed,middlearrow={triangle 45}] (1) to [bend right=15]  (3);

\end{tikzpicture}
  \caption{Orientations of lines in the  Fano plane induced by $\tau$ 
(red arrows).
  } 
\label{fig:OF}
\end{center}
\end{figure}
%%%%%%%%%%%%%%%%%%%%%
\noi
For example if $P\in\F$ and $\tau$ is of type $(0,1,3)$,  the orientation $(\tau^{\ast})^{-1}$ of $\F^\ast$  induces the  orientation on the line $L_P=\{D_P, D_{\tau^{-1}(P)}, D_{\tau^{-3}(P)} \}$ in $\F^\ast$
given by the cyclic permutation of lines in $\F$
\beqa
\label{eq:orD}
 (D_P ,D_{\tau^{-1}(P)}, D_{\tau^{-3}(P)} ) \ . 
\eeqa

\begin{proposition}\label{prop:orient}
 Let $(\F,\tau)$ and $(\F',\tau')$ be  oriented Fano planes, and let $f:\F\rightarrow\F'$ be an isomorphism of Fano planes. If $(\F,\tau)$ and $(\F',\tau')$ are isomorphic
(resp. not isomorphic)  then for all $D\in \F^\ast$, the restriction of $f $ to $D$ is  orientation preserving (resp. reversing) with respect to the induced orientations.
\end{proposition}

}

\section{Action of Aut$(\F)$ on composition factors and  the augmented Fano plane}\label{sec:cv-19}
In this section we first recall the definition and classification of composition factors which are the structures on $\F$
needed to define an eight-dimensional composition algebra \cite{mm}.
We then show that the action of $\A(\F)$ on composition factors has two orbits, each containing eight elements.
 These actions of $\A(\F)\cong GL(V_\F)$ are thus clearly  not isomorphic to the standard action of   $GL(V_\F)$ on
$V_\F$. We introduce a Radon transform taking functions on $\F$ to functions on $\F^\ast$. We show that its image (essentially) consists of
functions measuring the extent to which  elements of Aut$(\F)$ preserves  line orientations induced by composition factors.
Finally, we investigate the problem of lifting  automorphisms of $\F$  to
automorphisms of the associated composition algebra~\cite{wil}.

\subsection{Norms and multiplication factors}\label{sec:oct}

\begin{definition}
\begin{enumerate}\label{def:norm}
\item A norm on $\F$ is a function $N: \F \to \{-1,1\}$  such that
\beqa
N(P+ Q) = N(P) N(Q) \ , \ \ P \ne Q \ .\nn
\eeqa
\item Let $\F^2_0 = \{(P,Q) \in \F^2 \ \ \text{s.t} \ \ P \ne Q \}$. 
A multiplication factor is a map $\epsilon: \F^2_0 \to \{-1,1\}$  such that
$\epsilon_{PQ}+\epsilon_{QP}=0$. For $P \in \F$ the future  (resp. past) $\overrightarrow{P_\epsilon}$ (resp.
 $\overleftarrow{P_\epsilon}$)
of $P$ is defined by:
$\overrightarrow{P_\epsilon}=\{Q\in \F \ \text{s.t.} \ \ \epsilon_{PQ}=1\}$  (resp.
$\overleftarrow{P_\epsilon}=\{Q\in \F \ \text{s.t.} \ \ \epsilon_{PQ}=-1\}$). 
\item  Let ${\mathbb O_\F}$ be the set of $\FF-$valued functions  on the Fano cube $V_\F$.
\end{enumerate}
\end{definition}

\begin{remark}
It is easy to see that if $N$ is  a norm  then the set $\big\{P \in \F \ \mbox{s.t.}\  N(P)=1\big\}$ is either $\F$ or a line.
\end{remark}

For $P\in V_\F$  define $e_P\in {\mathbb O_\F}$ by
\beqa
e_P(Q) = \left\{
\begin{array}{ll}
1& \mathrm{~if~} Q=P\\
0& \mathrm{~if~} Q\ne P\ . 
\end{array}\right.\nn
\eeqa

Then $\{e_P \ \text{s.t.} \ P\in V_\F\}$ is an $\FF-$basis of ${\mathbb O_\F}$ and
\beqa
{\mathbb O_\F} = \FF e_0 \oplus \mathrm{Vect}\Big<e_P\ \text{s.t.}\  P \in \F\Big> \ . \nn
\eeqa
A norm  and a multiplication factor on $\F$ allow us to endow ${\mathbb O_\F} $ with a norm and a  multiplication:

\begin{definition}\label{def:mo}
%{\color{white} toto} 
Let $\F$ be a Fano plane equipped with a norm $N$ and a multiplication factor $\epsilon$.
\begin{enumerate}
\item
The multiplication $\cdot_\e : {\mathbb O_\F} \times {\mathbb O_\F} \to {\mathbb O_\F}$ is the unique bilinear map such that
\begin{enumerate}
\item For all $P\ne Q \in \F$:\ \ \ \  $e_P\cdot_\e e_Q = \epsilon_{PQ} e_{P + Q}$;
\item For all $P \in \F$:\ \ \ \   $e_P\cdot_\e e_P = -N(P)e_0$;
\item For all $P \in V_\F  $:\ \ \ \   $e_0 \cdot_\e e_P = e_P \cdot_\e e_0 = e_P$
\ \ (and we henceforth   denote $e_0$ by $1$).
\end{enumerate}
\item The norm $N_{\mathbb O_\F}: {\mathbb O_\F} \to \mathbb F$ is the  quadratic form:
$N_{\mathbb O_\F} (\lambda^0 e_0 + \sum \limits_{P\in\F} \lambda^P e_{P})= (\lambda^0)^2 + \sum\limits_{P\in \F} (\lambda^P)^2 N(P)$.

\end{enumerate}
\noindent We denote by ${\bf 1}$ the quadratic form on ${\mathbb O_\F}$ associated to the trivial norm on $\F$. By definition
 the triple  $({\mathbb O_\F},N_{\mathbb O_\F},\e )$ is a composition algebra  iff $N_{\mathbb O_\F} (Z  \cdot_\e W) =N_{\mathbb O_\F}(Z) N_{\mathbb O_\F}(W), \forall Z, W \in {\mathbb O_\F}$ and in this case we say that $\epsilon$ is a composition factor.
\end{definition}

The following proposition gives a necessary and sufficient condition for  $({\mathbb O_\F}, N_{\mathbb O_\F}, \epsilon)$ to be  a composition algebra (see \cite{mm}).

\begin{proposition}\label{prop:comp}
With the notation above $({\mathbb O_\F}, N_{\mathbb O_\F}, \epsilon)$ is a composition algebra  {\it iff}:
\beqa
%\label{eq:comp1}
(i)&&N(P+R) \epsilon_{PQ} \epsilon_{QR} = 1\hskip 4.1truecm \text{for any  line }\{P,Q,R\}, \nn\\
%\label{eq:comp2}
(ii)&& N(P+Q)\;\epsilon_{PQ}\epsilon_{QR}\epsilon_{RS}\epsilon_{SP} \;N(P+S) =-1\quad \text{for any  quadrilateral }\{P,Q,R,S\}.\nn
\eeqa
\end{proposition}

\begin{corollary}\label{cor:compt-or}
If $(\mathbb O_\F, {\bf 1}, \epsilon)$ is a composition algebra each line $\{P,Q,R=P+Q\}$ of $\F$ is canonically oriented by:
$\epsilon_{PQ}=\epsilon_{QR} =\epsilon_{RP}=1$.
\end{corollary}

%\begin{definition}
%Let $\F$ be Fano plane and let $N$ be a norm. 
%\begin{enumerate}
%\item {\color{red} garder est-ce qu'on l'utilise ???}
%If  $\epsilon$ is a composition factor,
%an $n-$gon $\{P_1,\cdots,P_n\}$  is $\epsilon-$orientable iff there exists a permutation $\sigma\in S_n$ such that $\epsilon_{P_{\sigma(1) }P_{\sigma(2)}} = \cdots=  \epsilon_{P_{{\sigma(n-1)}} P_{\sigma(n)}}= \epsilon_{P_{\sigma(n)} P_{\sigma(1)}}$.
%\end{enumerate}

%\end{definition}

\begin{example}\label{ex:com}
An oriented Fano plane $({\cal F},\tau)$ has a canonical composition factor $\epsilon^\tau$  for the norm ${\bf 1}$ defined using the Legendre symbol as follows. Fix 
$P_0\in \F$  and for $R,S\in\F$ set
\beqa
\epsilon_{RS}^\tau\ \  = \ \ \Big( \frac {j-i} 7 \Big)  \ \ = \ \  (j-i)^3 \ \  (\mathrm{mod} \  7),
\eeqa
where $R=\tau^i(P_0)$ and $S=\tau^j(P_0)$. This is well-defined,  independent of the choice of $P_0$ and
$\epsilon^\tau = \epsilon^{\tau^k}$ (resp. $\epsilon^\tau =- \epsilon^{\tau^k}$) if $k=1,2,4$ (resp. $k=3,5,6$).
Geometrically (see  Figure \ref{fig:OF})  this is equivalent to:
\beqa
\epsilon_{ P_iP_j}^\tau\ \  
%= \ \  (j-i)^3 \ \  (\mathrm{mod} \  7)\ \ 
= \ \ \left\{\begin{array}{ll}
\phantom{-}1 & \mathrm{if~there~is~an~arrow~from~} P_i \mathrm{~to~}  P_j\\
-1& \mathrm{if~there~is~an~arrow~from~} P_j \mathrm{~to~}  P_i  \nn
\end{array}
\right.  \ . 
\eeqa
One sees that every line is $\epsilon^\tau-$orientable and that there are  seven  $\epsilon^\tau-$orientable triangles $\{P_i, P_{i+2},P_{i+3}\}\, (i=1,\cdots,7)$.
Notice  that the complement of any line contains exactly one orientable triangle.
In fact these geometric properties are true for any composition factor for the norm {\bf 1} as shown in
\cite{mm}. Notice also that the orientation induced on each  line by $\tau$ \
(see Section \ref{sec:orientation}) is the same as the orientation on each line by $\epsilon^\tau$ (see Corollary \ref{cor:compt-or}).

%%%%%%%%%%%%%%%%%%%%%%%%%%%%%%%%%%%%%%%%%%%%%%%%%%%%%%%%%%%%%%
The multiplication table corresponding to $\epsilon^\tau$ 
 is then given  in Table \ref{tab:octmult}.
%%%%%%%%%%%%%%%%%%bbbbbbbbb
\begin{table}[!h]
          \caption{
          Octonion  multiplication and the oriented Fano plane.}
\begin{center}
%\begin{tabular}{@{}|c||c|c|c|c|c|c|c|c|@{}}
\begin{tabular}{cc}

\begin{minipage}{5.5cm}
\hskip .5truecm
%%%%%%%%%%%%%%%%%
\begin{tikzpicture}[scale=.7]

\tikzstyle{point}=[circle,draw]

\tikzstyle{ligne}=[thick]
\tikzstyle{pointille}=[thick,dotted]

\tikzset{->-/.style={decoration={
  markings,
  mark=at position .5 with {\arrow{>}}},postaction={decorate}}}

\tikzset{middlearrow/.style={
        decoration={markings,
            mark= at position 0.5 with {\arrow{#1}} ,
        },
        postaction={decorate}
    }}

\coordinate  (3) at ( -2,-1.15);%   {$\bullet$}; %[point] {$P_3$};%(-a/2,-asqrt(3)/6)

\coordinate  (2) at ( 0,-1.15) ;% {$\bullet$};%[point] {$P_2$};%(0,-asqrt(3)/6)
\coordinate  (5) at ( 2,-1.15);% {$\bullet$};% [point] {$P_5$};%(a,-asqrt(3)/6)

\coordinate  (6) at (0,2.31);% {$\bullet$};% [point] {$P_6$};%(0,a(sqrt(3)/2-sqrt(3)/6)
\coordinate  (1) at (1, .58);% {$\bullet$};% [point] {$P_1$};%(a/4,a(sqrt(3)/4-sqrt(3)/6)

\coordinate  (4) at (-1., .58) ;% {$\bullet$}; %[point] {$P_4$};
\coordinate  (7) at (0,0) ;%{$\bullet$};% [point] {$P_7$};

\draw [middlearrow={triangle 45}] (1) to [bend left=65] (2);
\draw [  middlearrow={triangle 45}] (2) --(3);
\draw [ middlearrow={triangle 45}] (3) --(4);
\draw [  middlearrow={triangle 45}] (6) --(7);
\draw [ middlearrow={triangle 45}] (7) --(1);
%\draw[color=red,dashed,middlearrow={triangle 45}] (4) to [bend left=15]  (5);

%\draw[color=red,dashed,middlearrow={triangle 45}] (5) to [bend right=15]  (6);
% -> latex oriente la droite
\draw [middlearrow={triangle 45}] (6) --(1);
\draw [ middlearrow={ triangle 45}] (1) --(5);

\draw[middlearrow={triangle 45}]  (5) --(2);

\draw[middlearrow={triangle 45}]  (3) --(7);

\draw  [middlearrow={triangle 45}]  (5) --(7);
\draw [middlearrow={triangle 45}] (7) --(4);
\draw[middlearrow={triangle 45}]  (4) --(6);
\draw [middlearrow={triangle 45}]  (7) --(2);

\draw [middlearrow={triangle 45}] (2) to [bend left=65] (4);
%\draw [middlearrow={triangle 45}] (4) to [bend left=65] (1);
\draw [middlearrow={triangle 45}]  (4) to [bend left=65] (1);
\draw (1)  node[scale=0.6] {$\bullet$} ;
\draw (1)  node {$\bullet$} ;
\draw (1) node[above=4pt  ]{$P_1$} ;

\draw (2)  node{$\bullet$} ;
\draw (2) node[below ]{$P_2$} ;
\draw (3) node[below left=0pt]{$P_3$} ;
\draw (3)  node{$\bullet$} ;
\draw (4)  node{$\bullet$} ;
\draw (4) node[ above left =1pt]{$P_4$} ;
\draw (5)  node{$\bullet$} ;
\draw (5) node[ below right ]{$P_5$} ;
\draw (6)  node{$\bullet$} ;
\draw (6) node[ above]{$P_6$} ;
\draw (7)  node{$\bullet$} ;
\draw (7) node[ above=5pt]{$P_7$} ;

%\draw[dashed,middlearrow={triangle 45}] (4) to [bend left=15]  (5);

%\draw[dashed,middlearrow={triangle 45}] (5) to [bend right=15]  (6);

%\draw[dashed,middlearrow={triangle 45}] (6) to [bend right=15]  (3);
%\draw[dashed,middlearrow={triangle 45}] (3) to [bend right=15]  (5);
%\draw[dashed,middlearrow={triangle 45}] (2) to [bend right=15]  (6);
%\draw[dashed,middlearrow={triangle 45}] (1) to [bend right=15]  (3);

\end{tikzpicture}

%%%%%%%%%%%%%%%%%%%%%%%%%%%%%
\end{minipage}
%%%%%%%%%%%%%%%%%%%%%%%%%%%%%%%%%%%%%%%%%%%%%%%%%%%
&
\begin{tabular}{@{}|c||c|c|c|c|c|c|c|c|@{}}
\hline
$ \curvearrowright$\hskip -.3truecm {\tiny $\cdot$}&$1$&$e_{P_1}$&$e_{P_2}$&$e_{P_3}$&$e_{P_4}$&$e_{P_5}$&$e_{P_6}$&$e_{P_7}$ \\
\hline\hline
$1$&$1$&$e_{P_1}$&$e_{P_2}$&$e_{P_3}$&$e_{P_4}$&$e_{P_5}$&$e_{P_6}$&$e_{P_7}$ \\
\hline
$e_{P_1}$&$e_{P_1}$&$-1$&$e_{P_4}$&$e_{P_7}$&$-e_{P_2}$&$e_{P_6}$&$-e_{P_5}$&$-e_{P_3}$ \\
\hline
$e_{P_2}$&$e_{P_2}$&$-e_{P_4}$&$-1$&$e_{P_5}$&$e_{P_1}$&$-e_{P_3}$&$e_{P_7}$&$-e_{P_6}$ \\
\hline
$e_{P_3}$&$e_{P_3}$&$-e_{P_7}$&$-e_{P_5}$&$-1$&$e_{P_6}$&$e_{P_2}$&$-e_{P_4}$&$e_{P_1}$ \\
\hline
$e_{P_4}$&$e_{P_4}$&$e_{P_2}$&$-e_{P_1}$&$-e_{P_6}$&$-1$&$e_{P_7}$&$e_{P_3}$&$-e_{P_5}$ \\
\hline
$e_{P_5}$&$e_{P_5}$&$-e_{P_6}$&$e_{P_3}$&$-e_{P_2}$&$-e_{P_7}$&$-1$&$e_{P_1}$&$e_{P_4}$ \\
\hline
$e_{P_6}$&$e_{P_6}$&$e_{P_5}$&$-e_{P_7}$&$e_{P_4}$&$-e_{P_3}$&$-e_{P_1}$&$-1$&$e_{P_2}$ \\
\hline
$e_{P_7}$&$e_{P_7}$&$e_{P_3}$&$e_{P_6}$&$-e_{P_1}$&$e_{P_5}$&$-e_{P_4}$&$-e_{P_2}$&$-1$ \\
\hline
\end{tabular}
%\label{tab:oct}
\end{tabular}
\end{center}
\label{tab:octmult}
\end{table}

\end{example}

Later on we will need the following proposition whose proof is immediate.
\begin{proposition}\label{prop:oct}
{\color{white}toto}
\begin{enumerate}

%\item
%\beqa
%\mathbb O = \mathbb F \oplus {\mathrm{Vect}}\Big<P \in \F \Big> \ .\nn
%\eeqa
\item  For any $P \in \F$, $e_P$ generates a  two-dimensional composition subalgebra of $\mathbb O_\F$.
\item 
If $(P,Q,R)\in\F^3$ are distinct aligned points then Vect$\big<1,e_P,e_Q,e_R\big>$ is an associative subalgebra isomorphic to a quaternion subalgebra
(which we also denote $\mathbb H_D$ where $D=\{P,Q,R\}$.)

\item 
If $(P,Q,R)\in \F^3$ are three  non-aligned points then $e_P,e_Q, e_R$  generate $\mathbb O_\F$.

\end{enumerate}
\end{proposition}

In \cite{mm} it was shown that if  $({\mathbb O_\F},{\bf 1} , \epsilon)$ is a composition algebra then
either $\overleftarrow{P}$ is a line for all $P\ \in\F$  or $\overrightarrow{P}$ is a line for
all $P\ \in\F$. Denoting ${\mathbb O_\F}_1^+$  (resp.  ${\mathbb O_\F}_1^-$)
the set of all composition algebras $({\mathbb O_\F}, {\bf 1},\epsilon)$ such that 
$\overrightarrow{P}$ (resp. $\overleftarrow{P}$)  is a line for all $P\in\F$, it was shown that 
${\mathbb O_\F}_1^+$  (resp.  ${\mathbb O_\F}_1^-$) is an affine space for      ${\cal S}_0^2(V_\F^\ast)$,
 the space of bilinear forms $B$ on $V_\F$  satisfying $B(P, P ) = 0, \forall P\in\F$.

In fact, by the following lemma ${\mathbb O_\F}_1^+$  and ${\mathbb O_\F}_1^-$ are affine $V_\F-$spaces.
\begin{lemma}\label{lem:affine}
The vector spaces ${\cal S}_0^2(V_\F^\ast)$ and $V_\F$ are $GL(V_\F)$ equivariantly  isomorphic.
\end{lemma}
\begin{demo}
Let $\wedge: V_\F\times V_\F \to (V_\F)^\ast$ be the  unique bilinear form such that
\beqa
P\wedge Q=\left\{\begin{array}{ll}
0&\hskip .15truecm \mbox{if}\ \ P=Q\ ,\\[5pt]
\begin{minipage}{5.55cm} {the unique non-zero linear form which vanishes at} $P$
 \mbox{and} $Q$ \end{minipage} &\begin{array}{l}\\  \mbox{if}\ \ P\ne Q \ .  \end{array}
\end{array}\right.\nn
\eeqa
One checks that the map $V_\F \to {\cal S}_0^2(V_\F^\ast)$ given by
\beqa
P \mapsto \Big[(Q,R)\mapsto \big(Q\wedge R\big)(P) \Big]\nn 
\eeqa
is a $GL(V_\F)$ equivariant isomorphism.
\end{demo}

\subsection{The action of Aut$(\F)$ on composition factors for the norm ${\bf 1}$}

In the rest of the paper we will assume that  the norm $N={\bf 1}$.

\begin{definition}[see \cite{mm}]\label{def:Opm}
An oriented map is a map $\alpha : \F \to V_\F^\ast$ such that 
\begin{enumerate}
\item[(a)] For all $P\in \F,   \alpha_P(P)=1\ .$
\item[(b)] If $P\ne Q \in \F$ then: \; $\alpha_P(Q)+ \alpha_Q(P)=1\  .$
\end{enumerate}
We denote $\F_0$ the set of all oriented maps. 

\end{definition}

 It was  shown in \cite{mm} that
`exponentiation'  defines a bijection from $\F_0$ to  ${\mathbb O_\F}_1^+$ (resp.  ${\mathbb O_\F}_1^-$).
Recall  the exponential  map $e: \mathbb Z_2 \to S_2$ and its inverse $\ell : S_2 \to \mathbb Z_2$  are defined  by
\beqa
\label{eq:exp-log}
 e^0=1\ \  , \ \ e^1=-1 \ \ 
\text{and} \ \
\ell(1)=0 \ \ , \ \ \ell(-1)=1 \ .
\nn
\eeqa

\begin{definition}
Let $\F$ be a Fano plane, let $\alpha$ be an oriented map, and let 
$\epsilon$  be a composition factor  and $g\in$Aut$(\F)$.
\begin{enumerate}
\item The oriented map $g\cdot \alpha$ is defined by
\beqa
(g\cdot \alpha)_P(Q) = \alpha_{g^{-1}\cdot P}(g^{-1}\cdot Q) \ \ \forall P, Q\in \F \ . \nn
\eeqa
\item The composition factor $g\cdot \epsilon$ is defined by
\beqa
(g\cdot \epsilon)_{PQ}= \epsilon_{g^{-1}\cdot P g^{-1}\cdot Q}  \ , \ \ \forall P\ne Q \in \F \ \ . \nn
\eeqa
This gives  left actions of Aut$(\F)$ on $\F_0$, ${\mathbb O}_\F^\pm$ that commute with the exponentiation  of Theorem 4.9  \cite{mm}.
\end{enumerate}
\end{definition}

\begin{theorem}\label{prop:blabla}
The actions of Aut$(\F)$ on  $\F_0$,  ${\mathbb O_\F}_1^+$ and    ${\mathbb O_\F}_1^-$ are transitive.
\end{theorem}  
\begin{demo}
To prove  the theorem it is is enough to prove that the action of  Aut$(\F)$ on  $\F_0$ is transitive.
We begin by by proving a series of lemmas.

\begin{lemma}\label{lem:1}
Let $P_0 \in \F$ and let $D_{P_0} = \{L \in \F^\ast, \ \ \text{s.t.} \ \ P_0 \not \in L \}$. Then the map $E_{P_0}:
\F_0 \to D_{P_0}$ defined by
\beqa
E_{P_0}(\alpha)= \alpha_{P_0} \   \nn 
\eeqa
is two-to-one surjective.
\end{lemma}
\begin{demo}
Let $\alpha,\beta \in \F_0$ be such that $E_{P_0}(\alpha)= E_{P_0}(\beta)$. By definition this means that
\beqa
\label{eq:1}
\alpha_{P_0}(P)= \beta_{P_0}(P) \ \ \forall P \in \F \ . 
\eeqa
Since $\F_0$ is an affine $V_\F-$space (see Lemma  \ref{lem:affine}) there exists a unique $Q \in V_\F$ such that $\alpha= \beta +P_0$,
{\it i.e}, 
\beqa
\label{eq:2}
\alpha_{P_0}(P)= \beta_{P_0}(P) + P_0 \wedge P(Q) \  \ \forall P \in \F \ . 
\eeqa
By Eqs.[\ref{eq:1}-\ref{eq:2}] we have $  P_0 \wedge P(Q)=0$ for all $P\in\F$, and by Lemma \ref{lem:affine} this implies that either $Q=0$ or $Q=P_0$. However Card$(\F_0)=8$ and Card$(D_{P_0})=4$ so this proves the lemma.
\end{demo}

\begin{lemma}\label{lem:2}
Let $P_0\in \F$ and let $\alpha, \beta\in\F_0$ be such that $\alpha(P_0)=\beta(P_0)$. Then there exists $g \in \A(\F)$ such that
$\beta= g\cdot \alpha$.
\end{lemma}
\begin{demo}
Since $\alpha(P_0) = \beta(P_0)$ by Lemma \ref{lem:1} either $\beta=\alpha$ or $\beta= \alpha +P_0$. In the first case
$\beta = \text{Id}\cdot \alpha$  so we can suppose that $\beta=\alpha + P_0$.  Let  $L=\{P_0,Q,R\}$ be any line through $P_0$ such
$\alpha_{P_0}(Q)=0$.  Define $g \in \A(\F)$ by
\beqa
g \cdot P = \left\{
\begin{array}{cl}
P&\text{if}\  P\in L\\
P+Q&\text{if}\ P \not \in L \ .
\end{array}
\right.\nn
\eeqa
Then one checks that $g\in\A(\F)$ and that $\beta=\alpha +P_0= g\cdot \alpha.$
\end{demo}

We now prove the theorem. Let $\alpha, \beta \in\F_0$ and let $P_0\in\F$. If $\alpha(P_0)=\beta(P_0)$ there exists $g\in \A(\F)$
such $\beta=g\cdot \alpha$ by Lemma \ref{lem:2}.  If $\alpha(P_0)\not =\beta(P_0)$, since $\alpha(P_0)$ and $\beta(P_0)$
are two lines not containing $P_0$ there exists $g\in\A(\F)$ such that $g\cdot P_0=P_0$
and $g$ maps the line $\alpha(P_0)$ to the line 
$\beta(P_0)$. Then $\beta$ and $g\cdot \alpha$ are two elements of $\F_0$ having the same value at $P_0$.
Again by Lemma \ref{lem:2} there exists $h\in\A(\F)$ such that $\beta=h\cdot \alpha$ and this completes the proof of the theorem. 
\end{demo}

\begin{remark}
Let $(\F,\tau)$ be an oriented Fano plane and fed 
The spaces $V_\F$ and $\F_0$ are affine $V_\F-$spaces and  the action of $\A(\F)$ on both is affine.  However the action
on the former is not transitive whereas the action
on the latter is transitive. This implies that there are two non-conjugate embeddings of $\A(\F)\cong GL(V_\F)$ in the  affine group
of $V_\F$. This phenomenon does not occur over $\mathbb R$ or $\mathbb C$, only over certain finite fields in certain
dimensions.
\end{remark}

\begin{remark}\label{rem:Je}
Let $(\F,\tau)$ be an oriented Fano plane and let $\epsilon^\tau$ be the canonical composition factor (see Example \ref{ex:com}).
By Theorem \ref{prop:blabla} the isotropy group ${\cal I}_{\e^\tau}= \{g\in \mathrm{Aut}(\F) \ \mathrm{s.t.} \ \ g\cdot \e^\tau = \e^\tau \}$
has $168/8=21$ elements of which six are of of seven, fourteen are of order three and the identity. 
The subgroup $\mathbb Z_7(\e^\tau)$ generated by $\tau$ is the cyclic group of order seven and one can show that
${\cal I}_{\e^\tau}$ is the normaliser of  $\mathbb Z_7(\e^\tau)$ in Aut$(\F)$. There is an exact sequence
\beqa
1 \to \mathbb Z_7(\e^\tau) \to {\cal I}_{\e^\tau} \to \mathbb Z_3 \to 1 \ , \nn
\eeqa
and elements of order three fall into two classes: (1) $z \in {\cal I}_{\e^\tau}$ such that $z\cdot g \cdot z^{-1}=g^2, \forall g \in \mathbb Z_7(\e^\tau)$ and (2) 
$z \in {\cal I}_{\e^\tau}$ such that $z\cdot g \cdot z^{-1}=g^4, \forall g \in \mathbb Z_7(\e^\tau)$.
For example one can take
\beqa
z=\begin{pmatrix} 1&2&3&4&5&6&7\\
                  4&1&5&2&6&3&7
                  \end{pmatrix} \ , \ \
z^2=  \begin{pmatrix} 1&2&3&4&5&6&7\\
                      2&4&6&1&3&5&7
                 \end{pmatrix} \    , \ \
\tau=  \begin{pmatrix} 1&2&3&4&5&6&7\\
                       2&3&4&5&6&7&1
                  \end{pmatrix}  \nn                
\eeqa
which satisfy
\beqa
z \cdot \tau \cdot z^{-1} = \tau^4 \ . \nn
\eeqa
\end{remark}

%In our presentation of this  construction  we have
%given a detailed study of composition factors.

\subsection{Radon transform}

In this subsection we introduce the  Radon transform associated to the incidence diagram
of points and lines in the Fano plane.

\begin{definition}
Let $\F$ be a Fano plane. The incidence space ${\cal I}$ is defined by
\beqa
{\cal I} = \Big\{(P,D) \in \F \times \F^\ast \ \ \text{s.t} \ \  P\in D\Big\} \ .
\nn
\eeqa
We denote ${\cal S}(\F)$ (resp. ${\cal S}(\F^\ast)$)  the set of functions from
$\F$ to $\mathbb Z_2$ (resp. $\F^\ast$ to $\mathbb Z_2$) and 
set  ${\cal S}_0(\F)= \big\{f \in {\cal S}(\F) \ \ \text{s.t.} \ \ \sum_{P\in \F}\ \ f(P)=0 \big\}$.

\end{definition}
We have the incidence diagram
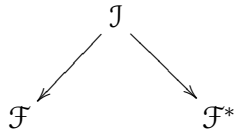
\beqa
\xymatrix{&{\cal I}\ar[dl]\ar[dr]\\
\F&&\F^\ast
}\nn
\eeqa
where the two maps are given by projections onto the first and second components respectively.
\begin{definition}
If $f\in {\cal S}(\F)$ we define  its  Radon transform  $f^\bigstar \in {\cal S}(\F^\ast)$  by
\beqa
f^\bigstar(D) =  \sum_{P\in D} f(P) \ . \nn
\eeqa
We set   $T(\F)= \big\{f \in {\cal S}(\F) \ \ \text{s.t.} \ \  f^\bigstar\equiv 0\big\}$.

\end{definition}

\begin{remark}
The Radon transform is a linear but not multiplicative map.
\end{remark}

The next proposition characterises the image and the kernel of the Radon transform.

\begin{proposition}\label{prop:kerim}
{\color{white}  toto}  \color{black}
\begin{enumerate}
\item  As a $\mathbb Z_2-$vector space  the kernel $T(\F)$ of the Radon transform is of dimension three.
\item  For each  $D\in \F^\ast$ define $T_D: \F \to \mathbb Z_2$ by
\beqa
T_D(P) = \left\{
\begin{array}{lcl}
0&\mathrm{if}&P\in D\\[4pt]
1&\mathrm{if}&P\not\in D
\end{array}\right. \ .
\nn
\eeqa
Then $T(\F)= \Big\{T_D \ \ \text{s.t} \ \  D\in \F^\ast\Big\}\cup \Big\{0\Big\}$.
\item 
Let $f \in {\cal S}(\F^\star)$. Then there exists $g\in {\cal S}(\F)$ s.t. $g^\bigstar=f$ iff
for any two triples $(D_1,D_2,D_3)$ and $(D'_1,D'_2,D'_3)$ of distinct concurrent lines,
we have $\sum_{i=1}^3 f(D_i)=\sum_{i=1}^3 f(D'_i)$.
\item The image $I$  of the Radon transform ${\scriptstyle  {}^{_\bigstar}} : {\cal S}(\F) \to {\cal S}(\F^\ast)$ is given by
\beqa
I = \Big\{T_P, T_P+1 \ \ \text{s.t.} \ \  \, P \in \F\Big\} \cup \Big\{0,1\Big\}  . \nn
\eeqa
Here $T_P$ is the function defined above corresponding to the line $P$ of the Fano plane $\F^\ast$.
\item Let $f \in {\cal S}(\F)$. There exists $P \in \F$ such that $f^\bigstar = T_P$ or $f^\bigstar \equiv 0$  {\it iff}
$f\in {\cal S}_0(\F)$.
In particular
\beqa
{\cal S}_0(\F)^\bigstar = \Big\{ f \in {\cal S}(\F^\ast) \ \ s.t. \ \ \sum^3_{i=1}  f(D_ i) = 0 \ \ {\mathrm{for~
  distinct~concurrent~lines~}} D_1,D_2,D_3\  \Big\} \ .
\nn
\eeqa
\end{enumerate}
\end{proposition}

\begin{demo}
1. Let $f$ be in $T(\F)$. Since its ´´average´´ on any line is zero it is easy to see that $f$ is completely determined 
by its values at three non-collinear
points (draw a picture).

2. Straightforward.

3. By the rank theorem the image of the Radon transform is of dimension 4 and cardinal 16. 
Let  $f$ be in ${\cal S}(\F)$. If $(D_1,D_2,D_3)$ are three distinct  concurrent lines, we have
\beqa
\label{eq:Lign}
\sum_{i=1}^3 f^\bigstar(D_i) %\sum_{P \in D_1 \cup D_2\cup D_3} f(P)
=3 f(D_1\cap D_2 \cap D_3) + \sum_{P \in \F \setminus  (D_1 \cap D_2\cap D_3)} f(P)
= \sum _{P\in \F} f(P)  
\eeqa
and therefore this sum is independent of the triple $(D_1,D_2,D_3)$. However this sum is either equal
to zero or  to one and there are only $8+8$ such functions in ${\cal S}(\F^\ast)$ by
 the argument of 1.

4. Straightforward.

5. This follows from \eqref{eq:Lign}.
\end{demo}

From now on we will use a multiplicative version of the Radon transform which we now explain. By abuse of language we refer to both
the additive and multiplicative versions as the Radon transform.
 In the statement of the following proposition we  use
the mutually inverse group isomorphisms  $e$ and $\ell$ given in \eqref{eq:exp-log}.

\begin{corollary}\label{lem:radon}
%Let ${\cal R}=\{h:\F \to S_2 \ \  \mathrm{s.t.} \  \prod_{P\in \F} h(P)=1\}$
Define ${\cal R}$ and ${\cal R}^\bigstar$ by:
\begin{enumerate}
\item[$\bullet$] ${\cal R}=e({\cal S}_0(\F))=\Big\{h:\F \to S_2 \ \  \mathrm{s.t.} \  \prod_{P\in \F} h(P)=1\Big\}$;
\item[$\bullet$] ${\cal R}^\bigstar= % =E(T(\F^\ast))$$\hskip .75truecm =
\Big\{h:\F^\ast\to S_2, \ \ \mathrm{s.t.}\ \ \prod_{i=1}^3 h(D_i)=1\ 
\mathrm{ for \ distinct  \ concurrent \ lines\ } \  D_1,D_2,D_3 \Big\}$.
\end{enumerate}
 Then $e\circ   {}^\bigstar\circ \ell: {\cal R} \to {\cal R}^\bigstar$
 is a group homomorphism and defines an exact sequence
\beqa
1 \to   e(T(\F)) \to {\cal R} \to {\cal R}^\bigstar \to 1 \ . \nn
\eeqa
In particular ${\cal R}$ is of  cardinal 64 and ${\cal R}^\bigstar$  of  cardinal 8.

\end{corollary}

\begin{demo}
After applying $e$, 
this follows from Proposition \ref{prop:kerim} (1) and (5). 
%By Proposition \ref{prop:kerim} we have an exact sequence of $\mathbb Z_2-$vector spaces:
%\beqa
%\xymatrix{
%0\ar[r]&T(\F) \ar[r]&{\cal S}(\F) \ar[r]^{\bigstar}&I\ar[r]& 0
%}\nn
%\eeqa
%which we can intersect with ${\cal S}_0(\F)$ to obtain the exact sequence
%\beqa
%\xymatrix{
%0\ar[r]&{\cal S}_0(\F)\cap T(\F) \ar[r]&{\cal S}_0(\F)\cap{\cal S}(\F) \ar[r]^{\hskip .8truecm \bigstar}&I'\ar[r]& 0 
%} \ . \nn
%\eeqa
%Since ${\cal S}_0(\F)\cap T(\F)=T(\F)$ it follows from Proposition \ref{prop:kerim} that $I'\subseteq I$ is of dimension $6-3=3$ and cardinal $8$. By \eqref{eq:Lign}
%we have $I' \subseteq L({\cal R}^\bigstar)$ and this inclusion is in fact an equality since the two sets have eight elements.
% Applying the isomorphism $E$  completes the proof of the proposition.
\end{demo}

\begin{remark}
The groups $e(T(\F))$ and ${\cal R}^\bigstar$ are  isomorphic to $(\mathbb Z_2)^3$.
\end{remark}
\noi
The eight  functions  of  ${\cal R}^\bigstar$  can be represented by the diagrams in Figure \ref{fig:delta*}.

%%%%% figure delta*
\input{delta-ast}
%label{fig:delta*}
%%%%%

%%%
\input{delta8}
%%%%

\begin{remark}
From Figure \ref{fig:delta*} we see that there is a natural one-to-one correspondence between ${\cal R}^\bigstar$ and $V_\F$, and
with respect to the natural $\mathbb Z_2-$vector space structures this is a vector space isomorphism. We have already seem that
$\F_0$ (see Definition \ref{def:Opm}) is affine space for $V_\F$ and so is also an affine space for ${\cal R}^\bigstar$. This action
is as follows: for $\alpha \in \F_0, f \in {\cal R}^\bigstar$ and $P\ne Q$, change $\alpha_P(Q)$ only if $f$ takes the value $-1$ on the line
through $P$ and $Q$.
\end{remark}
\begin{example}
If  $P\in\F$  denotes the point  at top of the triangle,
the fifth diagram in Figure \ref{fig:delta*} corresponds to the function $e(T_P)$.
The eight functions  in ${\cal R}$ whose Radon transform is $e(T_P)$ 
can be represented by the diagrams  in Figure \ref{fig:delta8}.
The point $P$ can be recovered from the diagrams representing these functions  as follows (see Figure \ref{fig:delta8}):
\begin{enumerate}
\item In  the second  diagram $P$ is the unique point  where the function is equal to $1$ (blue).
\item In the next four diagrams $P$ is the sum of the three blue points.
\item In the last three diagrams $P$ is the sum of the two red points.
\end{enumerate}
\end{example}

%%%%%%%%%%%%%%%%%%%%%%%%%%%%%%%%%%%%%%%%%%%%%%%%%%%%%%55
\subsection{The line orientation invariant $\delta^\ast$ of a composition factor}\label{sec:deltaast}

Let $\F$ be a Fano plane and $\e$ a composition factor.
In this section we introduce an invariant which measures whether or not an element of
Aut$(\F)$ preserves the  line orientations induced  by $\e$ (see \eqref{eq:orP}).

Let $g$ be an element of Aut$(\F)$ and let $D$ be a line
in $\F$.
If $P,Q$ are two distinct points of $D$ we temporarily set:
\beqa
\delta^\ast(g,P,Q) = \epsilon_{PQ} \epsilon_{g\cdot P g\cdot Q} \ , \nn
\eeqa
This number is either 1 or $-1$, and
we now show that it  is independent of the choice of $P,Q \in D$.

\begin{lemma}
With the notation above if $P',Q'$ are two distinct points of $D$, then
$\delta^\ast(g,P,Q) = \delta^\ast(g,P',Q')$,
\end{lemma}

\begin{demo}
Let $D=\{P,Q,R\}$. Without loss of generality we can suppose that $P'=Q$ and $Q'=R$.
From Proposition \ref{prop:comp} we have
\beqa
\epsilon_{PQ} \epsilon_{QR} = \epsilon_{g\cdot P g\cdot Q}  \epsilon_{g\cdot Q g\cdot R} =1 \ . \nn 
\eeqa
Since $\epsilon$ only takes the values $\pm 1$ we have
\beqa
\epsilon_{PQ}\epsilon_{g\cdot P g\cdot Q} = \epsilon_{QR} \epsilon_{g\cdot Q g\cdot R} \ , \nn
\eeqa
which proves the lemma
\end{demo}

This lemma shows that the following definition makes sense.
\begin{definition}\label{def:delta*}
Let $\F$ be a Fano plane and $\e$ a composition factor on $\F$.
Define $\delta^\ast : \mathrm{Aut}(\F) \times \F^\ast \to S_2$ by
\beqa
\delta^\ast(g,D) = \epsilon_{PQ} \epsilon_{g\cdot P g\cdot Q} \ , \nn
\eeqa
for any distinct points $(P,Q)$ of $D$.
\end{definition}

\begin{remark}\label{rem:mult}
It is immediate from the definition above that:
\beqa
\delta^\ast(g_2g_1,D)=\delta^\ast(g_2,g_1\cdot D) \delta^\ast(g_1,D) \ . \nn
\eeqa
\end{remark}

The algebraic interpretation of $\delta^\ast$ is as follows: if $g\in$Aut$(\F)$ and $D \in\F^\ast$ then $\delta^\ast(g,D)=1$ {\it iff} $g$ defines   an inclusion
of composition algebras $\mathbb H_D \hookrightarrow \mathbb O_\F$.
The geometrical interpretation  is:
 $\delta^\ast(g,D)=1$ {\it iff}  $g$ induces an orientation preserving map from $D$ to $g\cdot D$ (recall $\epsilon$ induces an orientation of $D$ and $g\cdot D$ see Eq.[\ref{eq:orD}]).

The function $(g,D)\mapsto \delta^\ast(g,D)$ is not arbitrary. In fact we will show that there are only eight possibilities
as a consequence  of the following two propositions.

\begin{proposition}\label{prop:det}
Let $\F$ be a Fano plane, $\e$ a composition factor on $\F$ and 
$g \in $ Aut$(\F)$. Then
$ \det g := \prod_{D \in \F^ \ast} \delta^\ast(g,D)=1$. 
\end{proposition}

\begin{demo}
Let $\F^\ast=\{D_1,\cdots,D_7\}$.
Let $g,h \in\text{Aut}(\F)$ and for each line $D_i$ choose distinct  $P_i,Q_i \in D_i$. Then
\beqa
\det g &= &\prod_{i=1}^7 \frac{\epsilon_{P_iQ_i}}{\epsilon_{g\cdot P_ig\cdot Q_i}} , \nn\\
\det h &= &\prod_{i=1}^7 \frac{\epsilon_{ P_i Q_i}}{\epsilon_{h\cdot P_ih\cdot Q_i}} . \nn
\eeqa
Choosing the point $(g\cdot P_i, g\cdot Q_i)$ on the line $g\cdot D_i$, the expression for $\det h$ can be written
\beqa
\det h &= &\prod_{i=1}^7 \frac{\epsilon_{g\cdot P_ig\cdot Q_i}}{\epsilon_{hg\cdot P_ihg\cdot Q_i}} . \nn
\eeqa
Thus
\beqa
\det g \det h = \det(hg)  \ ,\nn
\eeqa
and { det} defines a group homomorphism from Aut$(\F)$ to $S_2$. Since  Aut$(\F)$
is a finite simple group \cite{atlas},  this implies that
$\det g=1$ for any $g \in$Aut$(\F)$. 
\end{demo}

\begin{proposition}
Let $\F$ be a Fano plane, $\e$ be a composition factor on $\F$,
 $g\in$  Aut$(\F)$ and  $\{P,Q,R,S\}$ be a quadrilateral in $\F$.  Then
\beqa
\delta^\ast\big(g,P\wedge Q\big)\delta^\ast\big(g,Q \wedge R\big)= \delta^\ast\big(g,P\wedge S\big)\delta^\ast\big(g,S\wedge R\big) \ . \nn
\eeqa
\end{proposition}
\begin{demo}
By the quadrilateral rule (Proposition \ref{prop:comp}):
\beqa
\epsilon_{PQ} \epsilon_{QR} \epsilon_{RS} \epsilon_{SP} =
\epsilon_{g\cdot Pg\cdot Q} \epsilon_{g\cdot Qg\cdot R} \epsilon_{g\cdot Rg\cdot S} \epsilon_{g\cdot Sg\cdot P} =-1 \nn
\eeqa
which, using the definition of $\delta^\ast$, implies
\beqa
\delta^\ast\big(g,P\wedge Q\big)\delta^\ast\big(g,Q\wedge R\big)= \delta^\ast\big(g,P\wedge S\big)\delta^\ast\big(g,S\wedge R\big) \ .
\nn
\eeqa
\end{demo}
This proposition can be expressed in an equivalent more geometric  form.
For this let us remark the following property of lines in the Fano plane:
if $L$ is any line in the Fano plane there is a unique partition
$\{L_1,L_2\}, \{L_3,L_4\},\{L_5,L_6\}$
of the six remaining lines $L_1,\cdots, L_6$ such 
that  $L=\{L_1\cap L_2, L_3\cap L_4, L_5\cap L_6\}$.

\begin{proposition}\label{prop:delta*2}
Let $g \in$ Aut$(\F)$.
Let $L$ be a line of the Fano plane and let $L_1,\cdots,L_6$ be as above.
Then,
\begin{enumerate}
\item
$
\delta^\ast(g,L_1) \delta^\ast(g,L_2)=
\delta^\ast(g,L_3) \delta^\ast(g,L_4)=
\delta^\ast(g,L_5) \delta^\ast(g,L_6) 
$.
\item Let $P$ and $P'$ be two points of the Fano plane and  let
$D_1,D_2,D_3$ (resp.  $D'_1,D'_2,D'_3$) be the three lines passing through $P$
(resp. $P'$). Then
\beqa
\delta^\ast(g,D_1) \delta^\ast(g,D_2)\delta^\ast(g,D_3)=
\delta^\ast(g,D'_1) \delta^\ast(g,D'_2)\delta^\ast(g,D'_3) \ . \nn
\eeqa
\end{enumerate}
\end{proposition}

\begin{demo}
1: When the line $L$ is removed   from the Fano plane $\F$,  the six remaining lines are exactly the six sides
of the  quadrilateral $\F \setminus L$,  and the partition above $(L_1,L_2)$, $(L_3,L_4)$,  $(L_5,L_6)$ is obtained by grouping opposite sides of the quadrilateral.
We can certainly label the vertices of the quadrilateral  $\F \setminus L$ in such a way that:
$$L_1=P\wedge Q\ , \ \  L_2=R\wedge S\ , \ \  L_3=Q\wedge R\ , \ \ L_4=S\wedge P\ ,\ \  L_5=R\wedge P\ ,\ \  L_6=Q\wedge S \  $$
By the proposition above we have
\beqa
\delta^\ast(g,P \wedge Q)\delta^\ast(g,Q\wedge R)= \delta^\ast(g,P\wedge S)\delta^\ast(g,S\wedge R) \ , \nn
\eeqa
which implies
\beqa
\delta^\ast(g,L_1)\delta^\ast(g,L_3)= \delta^\ast(g,L_4)\delta^\ast(g,L_2) \ . \nn
\eeqa
Multiplying this equation by $\delta^\ast(g,L_3)\delta^\ast(g,L_2)$ gives
\beqa
\delta^\ast(g,L_1)\delta^\ast(g,L_2)= \delta^\ast(g,L_4)\delta^\ast(g,L_3) \ . \nn
\eeqa
This proves one of the  desired equalities, the others follow in a similar fashion.

\noi
 2: It is always possible to find two points $R$ and $S$ such that $\{P,P',R,S\}$ is a quadrilateral, {\it i.e.}, the complement of a line.
 The three lines passing through $P$ are $P\wedge P', P\wedge R,P\wedge S$ and the lines passing through $P'$ are
 $P\wedge P', P'\wedge R,P'\wedge S$. To prove 2 we have to show that
 \beqa
\delta^\ast(g,P\wedge P')\delta^\ast(g,P\wedge R)\delta^\ast(g,P\wedge S)=\delta^\ast(g,P\wedge P')
\delta^\ast(g,P'\wedge S)\delta^\ast(g,P'\wedge R) \ , \nn
\eeqa
which is equivalent to
\beqa
\delta^\ast(g,P\wedge R)\delta^\ast(g,P\wedge S) =\delta^\ast(g,P'\wedge S)\delta^\ast(g,P'\wedge R) \ , \nn
\eeqa
and this follows from 1.
\end{demo}

\begin{corollary}\label{cor:del}
Let $P \in \F$,  let $D_1,D_2,D_3$ be the three lines passing through $P$ and let $g \in$Aut$(\F)$. Then $\delta^\ast(g,D_1)\delta^\ast(g,D_2)\delta^\ast(g,D_3)=1$.
\end{corollary}
\begin{demo}
Suppose for contradiction that $\delta^\ast(g,D_1)\delta^\ast(g,D_2)\delta^\ast(g,D_3)=-1$. It follows from  Proposition \ref{prop:delta*2} 2 that the corresponding product  $\eta_{P'}$ for any point $P'$ in $\F$ is also equal to
$-1$ and hence that $\prod_{P' \in \F} \eta_{P'}=-1$. Since this product is also equal to det$(g)^3$ this contradicts Proposition \ref{prop:det} and the corollary is proved.
\end{demo}

\begin{theorem}\label{theo:delta*}
Let $(\F,\tau)$ be an oriented Fano plane, let $\e^\tau$ be the canonical composition factor
 and let $g\in$       Aut$(\F)$. Then $\delta^\ast(g,\phantom{D}): \F^\ast \to S_2$ is given by one of the eight diagrams
in Figure \ref{fig:delta*}. %In these diagrams $\delta^\ast$ takes the  value one on blue lines and $-1$ on red lines. In particular
%there are exactly eight different possibilities.
\end{theorem}
%%%%% figure delta*
%%%%%\input{delta-ast}
%label{fig:delta*}
%%%%%
\begin{demo}
By Corollary \ref{cor:del} the function $\delta^\ast(g,\dot) : \F^\ast \to S_2$ is in ${\cal R}^\bigstar$. By Proposition \ref{lem:radon} there are exactly eight possibilities, which
are the functions 
$e(T_P)$ for  $P \in \F$, and the constant function equal to one.
\end{demo}

%\begin{corollary}\label{cor:pd}
%The functions $\delta^\ast(g,\cdot)$ are characterised by the property: for any $P\in\F$ if $D_1,D_2,D_3$ are the three lines passing through $P$ then
%\beqa
%\prod_{i=1}^3\delta^\ast(g,D_i) = 1 \ . \nn
%\eeqa
%\end{corollary}

%\begin{demo}
%There are only eight functions on $\F^\ast$ satisfying this property and there are  eight functions $\delta^\ast(g,\cdot)$.
%\end{demo}

%\begin{remark}
%The   theorem    are true {\it mutatis mutandis} for any Fano plane and composition factor (see Theorem 
%\ref{prop:blabla}). 
%\end{remark}
%\begin{remark}\label{rem:point}
%There is a group structure on the set of eight diagrams in Figure \ref{fig:delta*} which can be seen in two ways. The first is to associate to each non-constant diagram
%the unique point of $\F$ which is the intersection of its blue lines. The constant diagram is then the identity, two distinct non-constant diagrams are composed according to the law
%of composition of the corresponding two points and the composition of a diagram with itself is the identity.

%The other way of describing this group structure is perhaps more surprising. From this point of view  we think of each diagram in Figure   \ref{fig:delta*} as representing a function
%on $\F^\ast$ which takes the value 0 on blue lines and 1  on red lines. Then one checks that the above group structure is  given by  point-wise addition of functions modulo two. 
%\end{remark}
\begin{example}
For the generators $a,b$ in \eqref{eq:ab} the functions $\delta^\ast(a,\cdot),\delta^\ast(b,\cdot)$ for the canonical composition factor $\e^\tau$ are given in Figure \ref{fig:ab}. 
%%%%%%%%%%%%%%%%%
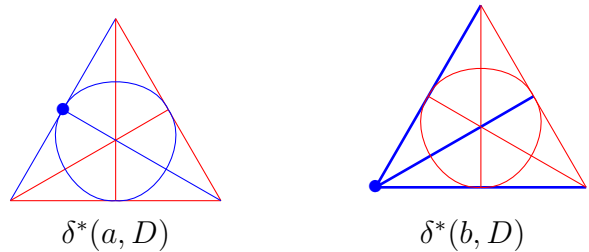
\begin{figure}[!ht]
\begin{center}
\begin{tabular}{ccc}
%%%%% 1
\begin{tikzpicture}[scale=.7]
\tikzstyle{point}=[circle,draw]

\tikzstyle{ligne}=[thick]
\tikzstyle{pointille}=[thick,dotted]

\tikzset{->-/.style={decoration={
  markings,
  mark=at position .5 with {\arrow{>}}},postaction={decorate}}}

\tikzset{middlearrow/.style={
        decoration={markings,
            mark= at position 0.5 with {\arrow{#1}} ,
        },
        postaction={decorate}
    }}

\coordinate  (3) at ( -2,-1.15);%   {$\bullet$}; %[point] {$e_3$};%(-a/2,-asqrt(3)/6)

\coordinate  (2) at ( 0,-1.15) ;% {$\bullet$};%[point] {$e_2$};%(0,-asqrt(3)/6)
\coordinate  (5) at ( 2,-1.15);% {$\bullet$};% [point] {$e_5$};%(a,-asqrt(3)/6)

\coordinate  (6) at (0,2.31);% {$\bullet$};% [point] {$e_6$};%(0,a(sqrt(3)/2-sqrt(3)/6)
\coordinate  (1) at (1, .58);% {$\bullet$};% [point] {$e_1$};%(a/4,a(sqrt(3)/4-sqrt(3)/6)

\coordinate  (4) at (-1., .58) ;% {$\bullet$}; %[point] {$e_4$};
\coordinate  (7) at (0,0) ;%{$\bullet$};% [point] {$e_7$};

% -> latex oriente la droite
\draw [color=red]  (6) --(1);
\draw [color=red] (1) --(5);

\draw [color=red]  (5) --(2);
\draw [color=red]  (2) --(3);

\draw [color=blue]  (3) --(4);
\draw  [color=blue] (4) --(6);

\draw  [color=red] (6) --(7);
\draw  [color=red] (7) --(2);

\draw [color=red] (3) --(7);
\draw [color=red]  (7) --(1);

\draw [color=blue] (5) --(7);
\draw [color=blue]  (7) --(4);

\draw [color=blue] (1) to [bend left=65] (2);
\draw [color=blue]  (2) to [bend left=65] (4);
%\draw [middlearrow={triangle 45}] (4) to [bend left=65] (1);
\draw  [color=blue] (4) to [bend left=65] (1);
%\draw (1)  node[scale=0.6] {$\bullet$} ;

\draw (4)  [blue] node{$\bullet$} ;

\end{tikzpicture}
&    \phantom{tttttt} &%%%%% 2
\begin{tikzpicture}[scale=.7]
\tikzstyle{point}=[circle,draw]

\tikzstyle{ligne}=[thick]
\tikzstyle{pointille}=[thick,dotted]

\tikzset{->-/.style={decoration={
  markings,
  mark=at position .5 with {\arrow{>}}},postaction={decorate}}}

\tikzset{middlearrow/.style={
        decoration={markings,
            mark= at position 0.5 with {\arrow{#1}} ,
        },
        postaction={decorate}
    }}

\coordinate  (3) at ( -2,-1.15);%   {$\bullet$}; %[point] {$e_3$};%(-a/2,-asqrt(3)/6)

\coordinate  (2) at ( 0,-1.15) ;% {$\bullet$};%[point] {$e_2$};%(0,-asqrt(3)/6)
\coordinate  (5) at ( 2,-1.15);% {$\bullet$};% [point] {$e_5$};%(a,-asqrt(3)/6)

\coordinate  (6) at (0,2.31);% {$\bullet$};% [point] {$e_6$};%(0,a(sqrt(3)/2-sqrt(3)/6)
\coordinate  (1) at (1, .58);% {$\bullet$};% [point] {$e_1$};%(a/4,a(sqrt(3)/4-sqrt(3)/6)

\coordinate  (4) at (-1., .58) ;% {$\bullet$}; %[point] {$e_4$};
\coordinate  (7) at (0,0) ;%{$\bullet$};% [point] {$e_7$};

% -> latex oriente la droite
\draw [color=red]  (6) --(1);
\draw[color=red]  (1) --(5);

\draw [line width=1pt,color=blue] (5) --(2);
\draw [line width=1pt, color=blue] (2) --(3);

\draw  [line width=1pt,color=blue]  (3) --(4);
\draw  [line width=1pt,color=blue]  (4) --(6);

\draw [color=red]  (6) --(7);
\draw  [color=red] (7) --(2);

\draw  [line width=1pt,color=blue]  (3) --(7);
\draw   [line width=1pt,color=blue]  (7) --(1);

\draw [color=red] (5) --(7);
\draw  [color=red] (7) --(4);

\draw [color=red] (1) to [bend left=65] (2);
\draw  [color=red] (2) to [bend left=65] (4);
%\draw [middlearrow={triangle 45}] (4) to [bend left=65] (1);
\draw [color=red]  (4) to [bend left=65] (1);
%\draw (1)  node[scale=0.6] {$\bullet$} ;

\draw (3)  [blue] node{$\bullet$} ;

\end{tikzpicture}

  \\

$\delta^\ast(a,D)$&&$\delta^\ast(b,D)$
\end{tabular}
  \end{center}

\caption{Values  of $\delta^\ast$ for the generators $a,b$ (blue lines:
  $\delta^\ast=1$,  red lines: $\delta^\ast=-1$).} 
\label{fig:ab}
\end{figure}
%%%%%%%%
\end{example}

%\begin{remark}
%From the explicit solutions above we see that the group Aut$(\F)$ has two orbits acting on the $\delta^\ast(g,\cdot)$'s\;: the constant solution and the non-constant solutions.
%\end{remark}

Recall that the group Aut$(\F)$ has 168 elements and that for fixed $g$ in Aut$(\F)$ there are 8 possibilities for the function $\delta^\ast(g,\cdot): \F^\ast \to S_2$ given by
the diagrams in Figure  \ref{fig:delta*}.
We now show that for each diagram there are exactly 21 elements of  Aut$(\F)$ corresponding to the same diagram.

\begin{proposition}\label{prop:delta*-glob}
Let $(\F,\tau)$ be an oriented Fano plane and let 
 $\epsilon^\tau$ be the canonical composition factor.  Let ${\cal I}_{\epsilon^\tau}= \big\{g\in \mbox{Aut}(\F): g\cdot \epsilon^\tau=\epsilon^\tau\big\}$ (see Remark \ref{rem:Je}).
\begin{enumerate}
\item If $g \in {\cal I}_{\epsilon^\tau}$ then $\delta^*(g,D)=1, \forall D \in \F^\ast$.
\item If $g \in {\cal I}_{\epsilon^\tau}$ and $h \in \mbox{Aut}(\F)$ then $\delta^\ast(gh,D)= \delta^\ast(h,D),$
and $\delta^\ast(hg,D)= \delta^\ast(h,g\cdot D),$
$\forall D \in \F^\ast$.
\end{enumerate}
\end{proposition}
\begin{demo}
1: This follows since
\beqa
\left\{
\begin{array}{l}
g\cdot \epsilon^\tau=\epsilon^\tau \ , \\
\delta^\ast(g,P\wedge Q) = \epsilon^\tau_{PQ} \epsilon^\tau_{g\cdot P g\cdot Q} \ .
\end{array}
\right. 
\nn
\eeqa

2: This follows from the  easily checked ``multiplier'' formula
\beqa
\label{eq:mult}
\delta^\ast(gh,D)=\delta^\ast(g,h\cdot D) \delta^\ast(h,D) \ ,
\eeqa
and 1 above.
\end{demo}

\begin{corollary}\label{cor:ttt}
For each of the 8 diagrams in Figure \ref{fig:delta*} there are 
 exactly 21 elements $g \in$ Aut$(\F)$ whose function $\delta^\ast(g,\cdot)$ corresponds to that diagram.
\end{corollary}

\begin{demo}
By Proposition \ref{prop:delta*-glob} (2) the map  $\delta^\ast$ factors to  a map from the left-coset space ${\cal I}_{\epsilon^\tau} \backslash \mbox{Aut}(\F)$ to the set of diagrams.
Note that the coset space has 8 elements since $|{\cal I}_{\epsilon^\tau}|=21$ (see Remark \ref{rem:Je}).
To show that this map  is a bijection  it is sufficient to show that it is  surjective.

By Proposition \ref{prop:delta*-glob},  we have $\delta^\ast(\tau,D)=1, \forall D \in \F^\ast$ and hence  $\delta^\ast(\tau,\cdot)$ corresponds to the first diagram in Figure \ref{fig:delta*}.
The other seven diagrams are characterised by the fact that their is a unique point $P$ such that $\delta^\ast(g,D)=1$ iff $P \in D$.
From this point of view the function $\delta^\ast(a,\cdot)$ corresponds to the point $P_4$ (See Figures \ref{fig:ab}).
However, by  \ref{prop:delta*-glob} (2),    the function    $\delta^\ast(a \tau^k,\cdot)$ corresponds to the  diagram with distinguished   point $\tau^{-k}(P_4)$,
and, since $\tau$ is of order seven, we obtain all seven diagrams in this way.
\end{demo}

\begin{remark}
By the multiplier  formula  above and  Eq.[\ref{eq:ab}] the function $\delta^\ast(\cdot,\cdot)$ is completely determined by  $\delta^\ast(a,\cdot)$ and $\delta^\ast(b,\cdot)$.
\end{remark}

We can now prove the converse of Proposition \ref{prop:orient}:
\begin{proposition}
Let $(\F,\tau)$ be an oriented Fano plane and let $\epsilon^\tau$ be the canonical composition factor.
If  $\tau'$ is an automorphism of order seven which induces the same orientation as $\tau$ on every line $D\in \F^\ast$,
then $\tau'=\tau,\tau^2$ or $\tau^4$.
\end{proposition}
\begin{demo}
If $\tau'$ induces the same orientation as $\tau$ on every line then $\delta^\ast(\tau',D) = \delta^\ast(\tau,D)=1$ for all $D\in\F^\ast$.
Since there are 21 elements of $\A(\F)$ with this  $\delta^\ast$ (see Corollary \ref{cor:ttt}) and since $|{\cal I}_{\epsilon^\tau}|=21$,
it follows from Proposition \ref{prop:delta*-glob} (1) that $\tau'\in {\cal I}_{\epsilon^\tau}$.
The only elements of  order seven in  ${\cal I}_{\e^\tau}$ are the six powers of $\tau$ (see Remark \ref{rem:Je})  and only $\tau'=\tau,\tau^2$ or $\tau^4$ induce the same  orientation on every line (Proposition \ref{prop:orient}).
\end{demo}

\subsection{Automorphisms of the augmented Fano plane}\label{sec:aug-auto}

%%%%%%%%%%%%%%%%%%%%%%%%%%%%%%%%%%%%%%%%%%%%%%%%%%%%%%%%%%%%%%%%%

In this section we investigate to what extent automorphisms of the Fano plane can be lifted to automorphisms of the octonions.

\begin{proposition}
Let $\F$ be a Fano plane and  let $\e$ be a composition factor.
Let $g\in \mathrm{Aut}(\F)$ and denote also by $g \in \mathrm{End}(\mathbb O_\F)$ the unique element such that
\beqa
g\cdot e_P= e_{g\cdot P} \ \ \forall P \in \F \ .  \nn
\eeqa
Let $D \in \F^\ast$ and let $\mathbb H_D$ denote the subalgebra of $\mathbb O_\F$ generated by the points of $D$ (see Proposition \ref{prop:oct} (2)). Then the following are equivalent:
\begin{enumerate}
\item $g|_{\mathbb H_D}$ is an algebra homomorphism;
\item $\delta^\ast(g,D)=1$.
\end{enumerate}
\end{proposition}
\begin{demo}
The proposition is a consequence of  the following equivalences: 
\beqa
&& g|_{\mathbb H_D} \ \  \mathrm{is ~an ~algebra ~homomorphism}\\
\Leftrightarrow&&g(e_P\cdot e_Q) = g(e_P)\cdot g(e_Q)\ ,   \forall P,Q \in D \nn\\
\Leftrightarrow&& \e_{PQ} e_{g\cdot(P+Q)}= \e_{g\cdot P g\cdot Q}e_{g \cdot P+g \cdot Q} \hskip .65truecm (cf. \ \ {\mathrm { Definition} \ \ \ref{def:mo}}) \nn\\
\Leftrightarrow&& \e_{PQ}=\e_{g\cdot P g\cdot Q}    \nn\\
\Leftrightarrow&&\delta^\ast(g,D)=1\hskip 3.55truecm   (cf. \ \ {\mathrm { Definition} \ \ \ref{def:delta*}})  \nn
\eeqa
\end{demo}
This means that  the canonical lift of $g \in$Aut$(\F)$ acts by automorphism on $\mathbb O_\F$ {\it iff} $\delta^\ast(g,D)=1, \forall D\in \F^\ast$.
By Corollary \ref{cor:ttt} there certainly exists $g\in$Aut$(\F)$ and $D\in\F^\ast$ such that $\delta^\ast(g,D)\ne 1$. However, we now exhibit a
finite group which acts on the octonions by automorphism and which is an eightfold non-trivial covering of Aut$(\F)$. To this end we introduce
the augmented Fano plane.

\begin{definition} \label{prop:wilson}
Let $\F= \{P_1,\cdots, P_7\}$ be a Fano plane, $\e$ be a composition factor and $\cdot$ be the associated composition product on $\mathbb O_\F$.
\begin{enumerate}
\item The  augmented Fano   plane
 $\hat \F_\epsilon$ is the subset of Im($\mathbb O_\F$) defined   by
\beqa
\hat \F_\epsilon =\Big\{\pm e_{P_1},\cdots, \pm e_{P_7}\Big\} \ . \nn
\eeqa
\item The group Aut$(\hat \F_\epsilon)$ of automorphisms of the augmented Fano plane is defined by
\beqa
\mbox{Aut}(\hat \F_\epsilon) = \Big\{ \hat g: \hat \F_\epsilon \to \hat \F_\epsilon \ \mbox{s.t.}\  &(i)&   \hat g(-e_P)=-\hat g(e_P)\  \    \forall P \in \F\ ;\nn\\
&(ii)& \hat g(e_P\cdot e_ Q) = \hat g(e_P)\cdot \hat g(e_Q) \ \ \forall P\neq  Q \in \F \Big\} \ . \nn 
\eeqa
\end{enumerate}
\end{definition}

\begin{remark}\label{rq:ImO}
Any  $\hat g:\hat \F_\epsilon \to \hat \F_\epsilon$ satisfying (i) above extends first  to an invertible  linear map $\tilde g:  \mathrm{Im}(\mathbb O_\F) \to \mathrm{Im}(\mathbb O_\F)$,
and then to  an invertible  linear map $g_{\mathbb O_\F}:\mathbb O_\F \to \mathbb O_\F$ by setting $g_{\mathbb O_\F}(1)=1$.
It is clear that $g \in \mathrm{Aut}(\hat \F_\epsilon)$ iff $g_{\mathbb O_\F}$ is an automorphism of $\mathbb O_\F$.
\end{remark}
There is a natural homomorphism  $\pi$ from Aut$(\hat \F_\epsilon)$ to Aut$(\F)$ such that
\beqa
\label{eq:Dcom}
\xymatrix{
\hat \F_\epsilon \ar[rr]^{\hat g} \ar[dd]^{/\mathbb Z_2}&& \hat \F_\epsilon \ar[dd]^{/\mathbb Z_2} \\\\
\F \ar[rr]^{\pi(\hat g)}&& \F
}  
\eeqa
is commutative and we have the following proposition: 
\begin{proposition}\label{theo:deltatilde}
Let $g \in $  {\rm{Aut}}$(\F)$ and let $\tilde \delta :  \F \to S_2$.
We define $\hat g:  \hat \F_\epsilon \to \hat \F_\epsilon$ by
\beqa
\hat g \cdot e_P = \tilde \delta(P) e_{g\cdot P} \ . \nn
\eeqa
Then $\hat g \in $ {\rm{Aut}}$(\hat \F_\epsilon)$ iff for any line $D=\{P,Q,R\}$ we have
\beqa
\label{eq:rad2}
\delta^\ast(g,D) = \tilde \delta(P)  \tilde \delta(Q)  \tilde \delta(R)  \ , 
\eeqa\
and in this case $\pi(\hat g) = g$.
\end{proposition}
\begin{demo}
Let $P, Q \in \F$ be two distinct points and let $R= P + Q$. Without loss of generality we can suppose that  $\epsilon_{PQ}=1$
and since $e_P \cdot e_Q=  e_R$  we have
\beqa
e_{g\cdot P} \cdot e_{g\cdot Q} = \epsilon_{rg\cdot P\, g\cdot Q}\;  e_{g\cdot R} =\delta^\ast(g,D) e_{g\cdot R}   \nn
\eeqa
where $D$ is the line $\{P,Q,R\}$.  It follows that
\beqa
\big(\hat g \cdot e_P\big) \cdot \big(\hat g \cdot e_Q\big)= \tilde \delta(P) \tilde \delta (Q) \; e_{g\cdot P} \cdot e_{g\cdot Q} =
\tilde \delta(P) \tilde \delta (Q)  \delta^\ast(g,D)  e_{g\cdot R} \ .\nn
\eeqa
On the other hand,
\beqa
\hat g \cdot \big(e_P\cdot   e_Q\big) =\hat g \cdot e_R =\tilde \delta(R) e_{g\cdot R} \ , \nn
\eeqa
which proves the result.
\end{demo}

\begin{corollary}\label{prop:tL}
Let $\F$ be a Fano plane and $\e$ be a composition factor.
\begin{enumerate}
\item 
For each $g \in$ {\rm Aut}$(\F)$ there exist exactly eight elements $\hat g \in$   {\rm Aut}$(\hat \F_\epsilon)$ such that $\pi(\hat g) = g$.

\item For $D \in \F^\ast$ we define $t_D: \hat \F_\epsilon \to \hat \F_\epsilon$ by (see also \cite{wil})
\beqa
t_D(\pm e_P) = \pm e\circ T_D(P) \; e_P= \left\{
\begin{array}{rl}
\pm e_P&\mathrm{if} \ P \in D \\
\mp e_P&\mathrm{if} \ P \not\in D \ .
\end{array}
\right.\nn
\eeqa
(Here  $e: \mathbb Z_2 \to S_2$ in $e \circ T_D$ denotes the exponential map.)
Then $t_D \in \mathrm{Aut}(\hat\F_\epsilon)$ and 
$\mathrm{Ker}(\pi) = \big\{\mathrm{Id}\big\} \cup \big\{    t_D \ \ \text{s.t.} \  \   D \in \F^\ast\big\}$.

\end{enumerate}
\end{corollary}
%We now calculate the  kernel of $\pi$  and show that  it is surjective \cite{wil}.

\begin{demo}
(1) Given $g \in${\rm Aut}$(\F)$, finding $\hat g \in${\rm Aut}$(\hat \F_\epsilon)$ such $\pi(\hat g) =g$  is equivalent to finding a function $\tilde \delta : \F
\to S_2$ such that Eq.[\ref{eq:rad2}] is satisfied. 
But in the language of  Corollary \ref{lem:radon} this means precisely   that the multiplicative Radon transform of
$\tilde \delta$ is $\delta^\ast(g,\cdot)$.  However, $\delta^\ast(g,\cdot) \in {\cal R}^\bigstar$ (see Corollary \ref{cor:del}) and again by
Corollary  \ref{lem:radon} any element of ${\cal R}^\bigstar$ is the
Radon transform of eight distinct elements of ${\cal R}$ (see Eq.[\ref{eq:Lign}]).

(2) We have  $\delta^\ast({\mathrm{Id}},D) =1$ for all $ D\in \F^\ast$ and we know that by Corollary \ref{lem:radon} $E \circ {}^\bigstar \circ L (f)=1$ {\it iff} $f \in e(T(\F))$.
This proves the result.
\end{demo}

\begin{remark}
 Recall from Example \ref{ex:ab} that the elements \
\beqa
a = \begin{pmatrix}1&2&3&4&5&6&7\\
         1&2&7&4&6&5&3
         \end{pmatrix} \ , \ \ 
b = \begin{pmatrix}1&2&3&4&5&6&7\\
                   2&7&4&6&5&3&1
         \end{pmatrix} \  \nn
\eeqa
generate Aut$(\F)$. One checks  (with the obvious abuse of notation) that
\beqa
\label{eq:abh}
\hat a = \begin{pmatrix}1&2&3&4&5&6&7\\
                   1&2&7&4&-6&5&-3
         \end{pmatrix} \ , \ \ 
\hat b = \begin{pmatrix}1&2&3&4&5&6&7\\
                   -2&7&4&6&5&3&-1
         \end{pmatrix} \  
\eeqa
are elements of Aut$(\hat \F_\epsilon)$ satisfying $\pi(\hat a)=a$ and  $\pi(\hat b)=b$.  
\end{remark}

\begin{remark}
A natural question is whether or not the exact sequence
\beqa
1 \rightarrow \text{Ker}(\pi) \rightarrow \mathrm{Aut}(\hat \F_\epsilon) \rightarrow   \mathrm{Aut}(\F) \rightarrow 1 \nn  
\eeqa
is split. In fact it is known  that it is not split \cite{Sand} and the easiest way to see this  as follows.
Consider the order four element $c \in \mathrm{Aut}(\F)$ given by
\beqa
c= \begin{pmatrix} 1&2&3&4&5&6&7\\
                   4&2&6&1&7&5&3
   \end{pmatrix} \ .\nn
\eeqa
One checks that
\beqa
\hat c= \begin{pmatrix} 1&2&3&4&5&6&7\\
                        -4&2&6&1&7&5&-3
   \end{pmatrix} \ ,\nn
\eeqa
is an order 8 element of $\mathrm{Aut}(\hat \F_\epsilon)$ satisfying $\pi(\hat c) = c$.
The other elements of $\pi^{-1}(c)$ are $t_{D} \circ \hat c$ where $D$ is any line of $\F$ (see Proposition \ref{prop:tL}) 
and, as one also checks, they too are of order 8.
This shows that the exact sequence above is not split.

Note that if $g \in \mathrm{Aut}(\F)$ is of order 2 (resp. 3) then $\pi^{-1}(g)$ consists of  four elements of order 2 and
four elements of order 4 (resp. four elements of order 3 and four elements of order 6). If  $g \in \mathrm{Aut}(\F)$ is of order 7  then $\pi^{-1}(g)$ consists  of elements of order 7.
\end{remark}

\section{The Lie algebra $\mathfrak g_2(\mathbb F)$}\label{sec:g2}
%%%%%%%%%%%%%%%%%%%%%%%%%%%%%%%%%%%%%%%%%%%%%%%%%%%%%%%%%%%%%
\par 
{\it
%In this section we fix an oriented  Fano plane $(\F,\tau)$ and the associated canonical composition factor   which will be denoted $\e$. 
 Throughout this section $\mathfrak{so}(n,\mathbb F)$ denotes the Lie algebra of $n\times n$  skew-symmetric matrices over $\mathbb F$. Recall that
$\mathfrak{so}(4,\mathbb F)$ is isomorphic to $\mathfrak{so}(3,\mathbb F)\times \mathfrak{so}(3,\mathbb F)$.
If $\FF$ admits a non-trivial quadratic extension $\hat \FF$  then  we denote by $\mathfrak{su}(n,\FF)$ the $\FF-$Lie algebra of traceless, anti-hermitian $n\times n$
matrices over $\hat \FF$.}\\

In this section we fix an oriented  Fano plane $(\F,\tau)$  {\it of type $(0,1,3)$} and the associated canonical composition factor   will be denoted $\e$.
It is well-known
that   when $\mathbb F = \mathbb R$  the group of automorphisms  of  $\mathbb O_\F$ is a simple, compact exceptional Lie group usually denoted $G_2$ (see for instance \cite{har}).
R. Wilson in  \cite{wil} gave an ``elementary'' construction of $\g_2$ (the Lie algebra of $G_2$) together with an action by automorphism  of a finite group of order $1344$. This is exactly the $\A(\hat \F_\epsilon)$ of Section \ref{sec:aug-auto}.
 In this section we shall express the commutation relations of $\g_2$
in terms of the incidence relations of the oriented Fano plane. 
As a consequence we will show that to each point $P$ of the Fano plane one
can associate a Cartan subalgebra  $\mathfrak h_P$ of  $\mathfrak g_2(\mathbb F)$. This association has the following remarkable properties which reflect the geometry of the Fano
plane $\F$: 
\begin{enumerate}[label=\alph*] 
\item[(a)] there is a  decomposition: $\mathfrak{g}_2(\mathbb F) = \bigoplus\limits_{P \in \F} \mathfrak h_P$ and 
 if $P\ne Q$ then $[\mathfrak h_P,\mathfrak h_Q]=\mathfrak h _{P+Q}$;
\item[(b)]  for each line $D$ in $\F$ we have $\mathfrak s_D =\bigoplus\limits_{P \in D} \mathfrak h_P$ is a Lie subalgebra of $\mathfrak g_2(\mathbb F)$ isomorphic to $\mathfrak{so}(3,\mathbb F) \times \mathfrak{so}(3,\mathbb F)$ and such that $\mathfrak h_P$ is a Cartan
subalgebra for  $P\in D$;
\item[(c)]   for each  point $P$ in $\F$ we have
\beqa
{\mathfrak s}_P=\big\{g \in \g_2(\FF) \ \ \text{s.t.}\ \  g\cdot P=0 \big\} \ , \nn
\eeqa
is a Lie subalgebra isomorphic to $\mathfrak{su}(3,\FF)$ or $\mathfrak{sl}(3,\FF)$ containing $\mathfrak h_P$ as a Cartan subalgebra;
\item[(d)]  if $P_1,P_2,P_3 \in \F$ are not aligned then $\mathfrak{h}_{P_1} \oplus \mathfrak{h}_{P_2} \oplus \mathfrak{h}_{P_3}$ generates $\mathfrak g_2(\mathbb F)$;

\item[(e)]  to the three lines containing a point $P$, one can associate  three elements of  of $\mathfrak h_P$ whose sum is zero and  which generate the root diagram of $\mathfrak{g}_2(\FF)$.
\end{enumerate}
If $\g$ is a simple Lie algebra, a decomposition $\g= \oplus_i \mathfrak h_i$ as an orthogonal sum of Cartan subalgebras satisfying $[\mathfrak h_i,\mathfrak h_j]\subseteq \mathfrak h_{k(ij)}$  is called a multiplicative orthogonal decomposition \cite{wil,KT}.\\

To prove these results
our main tool  is the observation that   $\mathbb O_\F$ is a space of spinors of Im$(\mathbb O_\F)$ with respect to octonion multiplication. 
This enables us to realise $\g_2$ inside a Clifford algebra, and use Clifford algebra techniques to perform calculations.

\subsection{Octonions, spinors of $\mathfrak{so}(7,\mathbb F)$ and
the Lie algebra $\mathfrak g_2(\mathbb F)$
}\label{sec:oct-spin}
Let $\mathbb O_\F$ be the  octonions defined by Table \ref{tab:octmult}, let Im$(\mathbb O_\F)$ be the purely imaginary octonions and let $B$ be the symmetric bilinear form associated to the norm ${\bf 1}$.
In this subsection we realise the spinor representation of $\mathfrak{so}(\text{Im}(\mathbb  O_\F),-B)$ in terms of octonion multiplication (see {\it e.g.} \cite{har} p. 121).

Define $\rho: \text{Im}(\mathbb  O_\F)\to \text{End}(\mathbb  O_\F)$ by
\beqa
\rho(x)(o)= x\cdot o\ ,  \ \ \forall x \in \text{Im}(\mathbb  O_\F)\ , \ \  \forall o \in \mathbb O_\F \ . \nn
\eeqa
Octonion multiplication is not associative but it is alternate ({\it i.e.}, the associator is antisymmetric) and this implies
\beqa
(\rho(x))^2= -B(x,x) \text{Id} \nn
\eeqa
and hence $\rho$ extends to a surjective  associative algebra homomorphism $\hat \rho$
with kernel ${\cal J}$ of dimension $2^6$:
\beqa
1 \to {\cal J} \to C(\mathrm{Im}(\mathbb O_\F),-B)  \buildrel{\hat \rho }\over\to \mathrm{End}(\mathbb O_\F) \to 1 \ . \nn
\eeqa
With respect to the well known decomposition
\beqa
 C(\text{Im}(\mathbb  O_\F),-B) = \bigoplus\limits_{k=0}^7 C^k \ .\nn
\eeqa
we have ${\cal J} \cap C^k=\{0\}$.
%of the Clifford algebra
%$\hat \rho: C(\text{Im}(\mathbb  O),-B)\to \text{End} (\mathbb O_\F)$ 
Furthermore,  $C^1$ is canonically identified with Im$(\mathbb O_\F)$ and 
$C^2$,  spanned by elements of the form $x y-yx \, (x,y \in C^1)$, is closed for the commutator and  isomorphic to $\mathfrak{so}(7,\mathbb F)$ as a Lie algebra.
The bracket of $C^2$  with $C^1$ defines an action of the Lie algebra $C^2$ on $C^1$ which is isomorphic to the vector representation
of $\mathfrak{so}(7,\mathbb F)$.

\begin{definition}\label{def:g2}
Set
\beqa
\mathfrak g_2(\mathbb F) = \big\{c \in C^2 \ \ \mathrm{s.t.} \ \ \hat \rho(c)(1)=0 \big\} \ . \nn
\eeqa
This is a Lie subalgebra of $\mathfrak{so}(7,\mathbb F)$.
\end{definition}

\begin{remark}
The Clifford algebra $C($Im$(\mathbb O_\F),-B)$ has a two-dimensional centre Vect$(1,\epsilon)$ where $\epsilon^2=1$. This implies that $1/2(1+\epsilon)\big( \mathfrak g_2(\FF)\big)$ is a Lie subalgebra
of $C($Im$(\mathbb O_\F),-B)$ isomorphic to $\mathfrak g_2(\FF)$. This embedding is used in \cite{toppan}.
\end{remark}

\begin{remark}
It can be shown   (see {\it e.g.} \cite{har} p. 122) that this Lie algebra is isomorphic to the Lie algebra of the group of automorphisms of $\mathbb O_\F$.
\end{remark}

%In the Fano plane (see Figure \ref{fig:OF}) observe that 
%\beqa
%e_{P_k} = \tau^{k-1}(e_{P_1}) \ , \  k=1,\cdots,7 \ .  \nn
%\eeqa
Recall that  ${\cal F} = \{P_i  \ , i \in \mathbb Z_7\}$ where  for all $i,j$ in $\mathbb Z_7$ we have $P_{i+j}=\tau^{j}(P_i)$  (see Figure \ref{fig:OF}).
The $e_{P_i}$ form    an orthonormal basis of $\text{Im}(\mathbb  O_\F)$ and
recall that  the action of $C^2$ on $C^1$ in this basis is given by
\beqa
\label{eq:so7}
\big[e_{P_iP_j},e_{P_k}\big] &=&\delta_{ik} e_{P_j} - \delta_{jk} e_{P_i}   \ ,\nn\\
\big[e_{P_iP_j},e_{P_kP_\ell}\big]&=&\delta_{ik} e_{P_jP_\ell} - \delta_{jk} e_{P_iP_\ell} +\delta_{i\ell} e_{P_k P_j} - \delta_{j\ell } e_{P_kP_i} \ .   
 \eeqa
 where $e_{P_i P_j} = \frac 1 2\; e_{{P_i}} e_{{P_j}}, i \ne j$.
 To simplify notation in what follows we set 
\beqa
 e_{P_i,P_j}= \hat \rho(e_{P_i P_j})= \frac14  \hat \rho[e_{P_i},e_{P_j}]=\frac 14\big( \hat \rho(e_{P_i}) \hat \rho(e_{P_j})- \hat \rho(e_{P_j}) \hat \rho(e_{P_i})\big), \ \ 1\le i\ne  j\le 7  \ . \nn
\eeqa
To each point $P_i$ of the Fano plane we now associate a Cartan subalgebra $\mathfrak{h}_{P_i}$  of $\mathfrak{g}_2(\mathbb F)$.  The idea is that each of the three lines containing  $P_i$ 
enables us to write $e_{P_i}$ as a product. For instance we have (see Table \ref{tab:octmult})
\beqa
e_{P_1}=e_{P_5}\cdot e_{P_6}=e_{P_3}\cdot e_{P_7} =e_{P_2}\cdot e_{P_4} \ . \nn 
\eeqa
Hence
\beqa
 \big(e_{P_3,P_7}-e_{P_5,P_6}\big)(1)= \big(e_{P_5,P_6} - e_{P_2,P_4}\big)(1) = \big(e_{P_2,P_4}- e_{P_3,P_7}\big) (1) =0 \ ,\nn
\eeqa
and $e_{P_3,P_7}-e_{P_5,P_6}, e_{P_5,P_6} - e_{P_2,P_4},e_{P_2,P_4}- e_{P_3,P_7}$ are in $\mathfrak{g}_2(\mathbb F)$ by Definition \ref{def:g2}.

The calculations above motivates the following definition in which we define a function $X: {\cal I} \to \g_2(\FF)$ where ${\cal I}=\{(P,D)\in \F \times \F^\ast \ \ \text{s.t.} \ \  P \in D \}$ is
the incidence space:
\begin{DP}\label{DP}

Given a point $P_i$ in $\F$ recall that the three lines containing $P_i$ (see Proposition \ref{cor:droites})   are:
\beqa
D_i=\{P_i,P_{i+1},P_{i+3}\}\ , \ \ %\nn\\
D_{i-1}= \{P_{i-1},P_i,P_{i+2}\}\ , \ %\nn\\
D_{i-3}= \{P_{i-3},P_{i-2},P_i\} \ . \nn
\eeqa
Define  $X_{P_i,D_i},X_{P_i,D_{i-1}},X_{P_i,D_{i-3}}\in \hat \rho(C^2)$ by
\beqa
\label{eq:X}
X_{P_i,D_i}&=&e_{P_{i+2},P_{i-1}}- e_{P_{i-3},P_{i-2}} \ , \nn\\
X_{P_i,D_{i-1}}&=& e_{P_{i-3},P_{i-2}}- e_{P_{i+1},P_{i+3}}\ , \\
X_{P_i,D_{i-3}}&=& e_{P_{i+1},P_{i+3}}- e_{P_{i+2},P_{i-1}} \ .\nn
\eeqa
Then
\begin{enumerate}
\item The action of $X_{P,D}$ on $e_Q$ is given by 
\beqa
\big[X_{P,D}, e_Q\big] = \left\{
\begin{array}{ll}
0& \mbox{if}\  Q  \in D\\
\epsilon_{P,Q} \; \epsilon^\ast_{ (PQ),D}  \; e_{P+Q} & \mbox{if}\  Q \not \in D \ ,
\end{array}\right.  \nn
\eeqa
where $\epsilon^\ast$ is the canonical composition factor of $(\F^\ast,\tau^{\ast \;-1})$ (see Proposition \ref{cor:droites}).
\item $X_{{P_i},D_i}+ X_{{P_i},D_{i-i}} + X_{{P_i},D_{i-3}}=0$;
\item  $X_{{P_i},D_i},X_{{P_i},D_{i-1}}$ and $X_{{P_i},D_{i-3}}$ are elements of $ \mathfrak{g}_2(\mathbb F)$;
\item[$\it 3'.$] $ \mathfrak{g}_2(\mathbb F)=$Vect$\Big<X_{P_i,D_i},X_{P_i,D_{i-1}},X_{P_i,D_{i-3}} : i\in \mathbb Z_7\Big>$;
\item[$\it 4.$] The subspace $\mathfrak{h}_{P_i}= \mathrm{Vect}\Big<X_{{P_i},D_i},X_{{P_i},D_{i-1}},X_{{P_i},D_{p-i}}\Big> $ is a Cartan subalgebra of $\mathfrak g_2(\FF)$.

\end{enumerate}
\end{DP}
\begin{demo}
Straightforward. 
\end{demo}

\begin{remark}
Using only the incidence geometry of $\F$ and the composition factor $\epsilon$, on the $\FF-$vector space Im$(\mathbb O_\F)$ 
we can define the three-form
\beqa
\omega= \sum \limits_{D\in \F^\ast} \epsilon_{PQ} \;\epsilon_{QR} \;\epsilon_{RP}\; e^P\wedge e^Q \wedge e^R \nn 
\eeqa
where $\{e^1,\cdots, e^7\}$ is the  basis dual to $\{e_1\cdots, e_7\}$ and $D=\{P, Q, R\}$.
Note that each term in the above sum is independent of the order of the points chosen  in the corresponding line.
This is the unique three-form invariant
under $G_2$,  normalised such that $B(\Omega, \Omega)=7$ and oriented such that 
\beqa
i_v \omega\wedge i_w \omega \wedge \omega = -6 B(v,w) \; e^1 \wedge \cdots \wedge e^7 \ . \nn
\eeqa
These are the  same conventions as \cite{Ka} but opposite to \cite{br}.
Dually, again using only the incidence geometry of $\F$ and the composition factor $\epsilon$, we can define the four-form
\beqa
\Omega = \sum \limits_{K\in {\cal Q}} \epsilon_{PQ} \;\epsilon_{RS} \; e^P\wedge e^Q \wedge e^R \wedge e^S \nn 
\eeqa
where ${\cal Q}$ is the set of (seven) quadrilaterals of $\F$  and $K=\{P,Q,R,S\}$.
Again each term in the above sum is independent of the order of the points  chosen   in the corresponding quadrilateral.
This is the unique $G_2-$invariant four-form satisfying
\beqa
\Omega\wedge \omega= - 7 e^1 \wedge \cdots \wedge e^7 \ . \nn
\eeqa
These are the  same conventions as  \cite{br} but opposite to \cite{Ka}.
\end{remark}
\begin{remark}
Not only the $X_{P,D}$ but also their action on Im$(\mathbb O_\F)$ have now been given in terms of the incidence geometry of $\F$ and the
composition factors $\epsilon$ and $\epsilon^\ast$.
\end{remark}

The three generators of the Cartan subalgebra associated to a point  of $\F$ (see  Definition/Proposition \ref{DP}) give rise to a root system of type $G_2$ in the following sense 
\begin{proposition}
Let $X_{P_i,D_i},X_{P_i,D_{i-1}},X_{P_i,D_{i-3}}$ and $\mathfrak{h}_{P_i}$ be as in \ref{DP} and set:
\beqa
\label{eq:Y}
\begin{array}{lllll}
Y_{P_i,D_i}&=&X_{{P_i},D_{i-1}}-X_{P_i,D_{i-3}}&=& e_{P_{i-3},P_{i-2}} + e_{P_{i+2},P_{i-1}}-2 e_{P_{i+1},P_{i+3}} \\
Y_{P_i,D_{i-1}}&=&X_{P_i,D_{i-3}}-X_{P_i,D_{i}}&=& e_{P_{i-3},P_{i-2}} + e_{P_{i+1},P_{i+3}}-2 e_{P_{i+2},P_{i-1}}\\
Y_{P_i,D_{i-3}}&=&X_{P_i,D_{i}}-X_{P_i,D_{i-1}}&=&   e_{P_{i+1},P_{i+3}}+ e_{P_{i+2},P_{i-1}}-2 e_{P_{i-3},P_{i-2}}\ . \\
\end{array}
\eeqa
Then
\beqa
{\cal W} = \big\{\pm X_{P_i,D_i},\pm X_{P_i,D_{i-1}},\pm X_{P_i,D_{i-3}}, \pm Y_{P_i,D_i},\pm Y_{P_i,D_{i-1}},\pm Y_{P_i,D_{i-3}}\big\} \subset \mathfrak{h}_{P_i} \nn 
\eeqa
is a root system of type $G_2$. For instance if 
$\alpha= X_{P_i,D_{i}}$ and $\beta=Y_{P_i,D_{i-1}}$ then
\beqa
{\cal W} = \big\{\pm \alpha,\pm \beta, \pm(\alpha + \beta), \pm(\beta+ 2 \alpha), \pm(\beta+3\alpha), \pm(2 \beta + 3 \alpha)\} \ . \nn
\eeqa
\end{proposition}

\begin{demo}
This follows from the fact that  $\{e_{ij},\ \  1\le i <j\le 7 \}$ is an orthonormal basis of $C^2$.
\end{demo}

\subsection{Action of Aut$(\hat \F_\epsilon)$ on the $X_{P,D}$}\label{sec:Auth}
In this section we analyse the action of the group  Aut$(\hat \F_\epsilon)$ (see Definition \ref{prop:wilson})  on the elements $X_{P,D}$ of $\g_2(\mathbb F)$ introduced above.
Recall first that the $X_{P.D}$ are  elements of End$(\mathbb O_\F)$  ({\it c.f.} Section \ref{sec:oct-spin}) and
 the elements of the group
Aut$(\hat \F_\epsilon)$  act naturally on this vector space ({\it c.f.}  Remark \ref{rq:ImO}).  This means that it makes sense to conjugate elements of $\hat \rho(C^2)$ by  Aut$(\hat \F_\epsilon)$.

We start by making the following observation:
\begin{lemma}
Let $P \in \F, D\in \F^\ast$ and let $\hat g \in \mathrm{Aut}(\hat \F_\epsilon)$ be such $\pi(\hat g) = g \in  \mathrm{Aut}(\F)$. Then  {\underline{up to a sign}}:
\beqa
\hat g \; X_{P,D} \;\hat g^{-1} =  X_{g\cdot P, g\cdot D} \ . \nn
\eeqa

\end{lemma}
\begin{demo}
Let $Q$ be in $\F$. Then up to a sign,  the  LHS acting on $e_Q$ is given by  (see Eq.\ref{eq:Dcom})
\beqa
\hat g \; X_{P,D} \;\hat g^{-1}(e_Q)= \left\{
\begin{array}{ll}
0& \text{if~} \ g^{-1} \cdot Q \in D \\
\hat g \cdot e_{g^{-1}\cdot Q +P}= e_{Q+g\cdot P}&  \text{if~} \ g^{-1} \cdot Q \not\in D \ . \\
\end{array}\right.\nn
\eeqa
On the other hand  up to a sign, the RHS acting on $e_Q$ is given by
\beqa
 X_{g\cdot P, g\cdot D}(e_Q)
 =\left\{
 \begin{array}{ll}
0& \text{if~} \  Q \in g\cdot D \\
 e_{Q+g\cdot P}&  \text{if~} \  Q \not\in g\cdot D \ . \\
\end{array}\right.\nn
\eeqa
Comparing these two expressions  proves the lemma.
\end{demo}

This lemma means  the following definition is legitimate:
\begin{definition}\label{def:delta}
Let $(P,D)  \in {\cal I}$ be an incident pair and let $\hat g \in \mathrm{Aut}(\hat \F_\epsilon)$ be such $\pi(\hat g) = g \in  \mathrm{Aut}(\F)$.
Define $\delta : \mathrm{Aut}(\hat \F_\epsilon) \times {\cal I}  \to \mathbb Z_2$ by
\beqa
\hat g \; X_{P,D} \;\hat g^{-1} = \delta(\hat g,P,D) X_{g\cdot P, g\cdot D} \ . \nn
\eeqa
\end{definition}

We now show that in fact $ \delta$ does not depend on the line $D$.

\begin{proposition}
Let $P \in \F$ let $ D, D' \in \F^\ast$ be such that $P \in D\cap D'$. Let $\hat g \in \mathrm{Aut}(\hat \F_\epsilon)$ be such $\pi(\hat g) = g \in  \mathrm{Aut}(\F)$. Then
$
 \delta(\hat g,P,D)= \delta(\hat g,P,D') \ . \nn
$
\end{proposition}

\begin{demo}
Suppose  $D\ne D'$.
By proposition \ref{DP} if $D''=D + D'$ we have
\beqa
X_{P,D} + X_{P,D'} + X_{P,D''} &=&0 \nn\\
X_{g\cdot P,g\cdot D} + X_{g\cdot P,g\cdot D'} + X_{g\cdot P,g\cdot D''} &=&0 \ . \nn
\eeqa
Conjugating the first equation above by $\hat g$ we obtain
\beqa
\delta(\hat g,P,D)X_{g\cdot P,g\cdot D} + \delta(\hat g,P,D')X_{g\cdot P,g\cdot D'} + \delta(\hat g,P,D'')X_{g\cdot P,g\cdot D''} &=&0 \nn
\eeqa
and comparing this with the second equation above it follows that $\delta(\hat g,P,D)=\delta(\hat g,P,D')=\delta(\hat g,P,D'')$.
\end{demo}
From now on we write $\delta(\hat g,P)$   for $\delta(\hat g, P, D)$.
\begin{example}\label{ex:tL}
For $D\in \F^\ast$ and $t_D \in \mbox{Aut}(\hat \F_\epsilon)$ (see Corollary \ref{prop:tL}) one checks easily that 
\beqa
\delta(t_D,P) = \left\{ 
\begin{array}{rl}
1&\mathrm{if} \ P \in D\\
-1&\mathrm{if} \ P \not\in D\\
\end{array}= e^{T_D(P)}
\right.
\nn
\eeqa
where $T_D$ is defined in Proposition \ref{prop:kerim}.
We denote by $T$ the set of  all $T_D$ together with the constant function equal to one.
Note that $T$ is a group for pointwise multiplication.
\end{example}

\noi
%The next proposition shows that
%the function  $\delta$ has a multiplier property:
%\begin{proposition}
%Let $P\in \F, \hat g,\hat h \in \mathrm{Aut}(\hat \F_\epsilon)$ then
%\beqa
%\delta(\hat g\hat h ,P)= \delta(\hat g, \hat h \cdot P)\delta(\hat h ,P) \ . \nn
%\eeqa
%\end{proposition}
%\begin{demo}
%Straightforward.
%\end{demo}

%We now come to one of our  main result of this section which establish
We now establish a connection between the line orientation $\delta^\ast$ (see Section \ref{sec:deltaast}) and the function $\delta$.

\begin{proposition}\label{prop:espdel}
Let $(P, D)\in  \F \times \F^\ast$ and $\hat g \in \mathrm{Aut}(\hat\F_\epsilon)$ be such that $\pi(\hat g)=g$. Then
\beqa
\label{eq:pe}
\delta^\ast(g,D) = \prod_{P\in D} \delta(\hat g,P) \ .
\eeqa
Hence $\delta^\ast(g,\cdot)$ is the (multiplicative) Radon transform of $\delta(\hat g,\cdot)$.
\end{proposition} 

\begin{demo}
The proposition  is a consequence of the following two lemmas and the fact that $\mathrm{Aut}(\hat \F_\epsilon)$ is generated by $\hat a, \hat b$
given in \eqref{eq:abh} 
and the $t_{D_i}$ introduced in Proposition \ref{prop:tL}.
\begin{lemma}
If the formula \eqref{eq:pe} is true for $\hat g$ and $\hat h$, then it is true for $\hat g \hat h$.
\end{lemma}

\begin{demo}
By the multiplier  property of $\delta^\ast$ (see Remark \ref{rem:mult}) we have
\beqa
\delta^\ast(gh,D) =\delta^\ast(g,h \cdot D) \delta^\ast(h, D) \nn 
\eeqa
which is equal to
\beqa
 \prod_{P\in D} \delta(\hat g, h \cdot P)  \prod_{P\in D} \delta(\hat h,P) \nn
\eeqa
since  formula \eqref{eq:pe} is true for $\hat g$ and $\hat h$.
By the analogous multiplier property of $ \delta$ this is in turn equal to
\beqa
\prod_{P\in D} \delta(\hat g\hat h ,  P) \ .\nn
\eeqa
\end{demo}
\begin{lemma}
The formula  \eqref{eq:pe} is true for
\beqa
\hat a = \begin{pmatrix}1&2&3&4&5&6&7\\
                        1&2&7&4&-6&5&-3
                        \end{pmatrix} \ , \ \ \  
\hat b=\begin{pmatrix}1&2&3&4&5&6&7\\
                        -2&7&4&6&5&3&-1
                        \end{pmatrix}\nn
\eeqa
and all translations $t_{D_i}$.
\end{lemma}
\begin{demo}
To calculate $\delta(\hat a, \cdot)$ and  $\delta(\hat b, \cdot)$ one has to use Definition \ref{def:delta} and unfortunately this is tedious if straightforward. The computation of $\delta^\ast( a, \cdot)$ and  $ \delta^\ast( b, \cdot)$ is much easier but
we omit the details of both calculations. The results are given in
 Figure \ref{fig:eps} from which the lemma follows. 
\input{eps}
\end{demo}

Formula \eqref{eq:pe} in the general case follows from the lemmas above.
\end{demo}

\begin{corollary}
For all $\hat g$ in Aut$(\hat\F_\epsilon)$ we have $\delta(g,\cdot) \in {\cal R}$, {\it i.e.}, 
\beqa
\prod_{P \in \F} \delta(\hat g, P)=1 \ .\nn
\eeqa
\end{corollary}
\begin{demo}
This follows from the proposition above and  equation \eqref{eq:Lign}.
\end{demo}

\begin{remark}
Given $\hat g \in  {\mathrm{Aut}}(\hat \F_\epsilon)$ such that $\pi(\hat g)=g$  we have seen that the Radon transform of $\delta(\hat g,\cdot)$ is $\delta^\ast(g,\cdot)$.
According to Theorem \ref{theo:deltatilde},
this means that that if we define $\hat g^{\color{white}{t}}{}^{\hskip -.1truecm '}: \hat \F_\epsilon \to \hat \F_\epsilon$  by
\beqa
\hat g^{\color{white}{t}}{}^{\hskip -.1truecm '}\cdot e_{P} = \delta(\hat g ,P) e_{g\cdot P}  \ , \nn
\eeqa
then $\hat g^{\color{white}{t}}{}^{\hskip -.1truecm '} \in {\rm{Aut}}(\hat \F_\epsilon)$ and
$\pi(\hat g^{\color{white}{t}}{}^{\hskip -.1truecm '})=g$. Thus by Corollary \ref{prop:tL} either
$\hat g'= \hat g$ or there exists $D\in \F^\ast$ such that
\beqa
\hat g^{\color{white}{t}}{}^{\hskip -.1truecm '}= t_D \circ \hat g \ . \nn
\eeqa
Both cases occur.
For example if $\hat g = t_D$ then by Example \ref{ex:tL}, we have $\hat g^{\color{white}{t}}{}^{\hskip -.1truecm '}= t_D$.
On the other hand  if
\beqa
\hat a &=& \begin{pmatrix}1&2&3&4&5&6&7\\
                        1&2&7&4&-6&5&-3
                        \end{pmatrix}\nn
\nn
\eeqa
then $\delta(\hat a,\cdot)$ is given by the first diagram in Figure \ref{fig:eps} and
therefore
\beqa
\hat a ^{\color{white}{t}}{}^{\hskip -.1truecm '}=
\begin{pmatrix}1&2&3&4&5&6&7\\
               1&-2&-7&-4&-6&5&3
                        \end{pmatrix}=
t_{D_5} \circ \hat a \ . \nn
\eeqa

\end{remark}

%%%%%%%%%%%%%%%%%%%%%%%%%%%%%%%%%%%%%%%%%%%%%%

%%%%%%%%%%%%%%%%%%%%%%%%%%%%%%%%%%%%%%%%%%%%%%%%%%%%%
\subsection{Commutation relations for $\g_2$}

The incidence properties of the oriented Fano plane define multiplication of the octonions (see  Table \ref{tab:octmult}).
In an analogous fashion they also define the Lie bracket of  $\g_2(\mathbb F)$.
%and the following proposition is the  first step in expressing   the bracket of $\g_2(\mathbb F)$  entirely in terms of the oriented Fano plane. 
In this section we will express all brackets of the $X_{P,D}$'s in $\g_2(\mathbb F)$ in terms of the incidence relations of the oriented Fano plane.
To begin we need the following lemma:

\begin{lemma}\label{lem:inci}
Let ${\cal I}$ be the incidence space
\beqa
{\cal I} = \Big\{ (P,D) \in \F \times \F^\ast \ \ \text{s.t.} \ \  P \in D \Big\} \ . \nn
\eeqa
Then the orbits of Aut$(\F)$ acting on ${\cal I} \times {\cal I}$ are given by
\beqa
{\cal I} \times {\cal I} = {\cal D} \cup {\cal O}_1 \cup  {\cal O}_2 \cup  {\cal O}'_3 \cup {\cal O}_3 \cup  {\cal O}_4 \ , \nn
\eeqa
where
\beqa
{\cal D} &=& \Bigg\{\Big((P,D),(P',D')\Big) \in {\cal I}\times {\cal I}: \  P=P',D=D' \Bigg\} \nn\\
{\cal O}_1 &=& \Bigg\{\Big((P,D),(P',D')\Big)\in {\cal I}\times {\cal I}:\   P=P', D \ne D'\Bigg\} \nn\\
{\cal O}_2 &=& \Bigg\{\Big((P,D),(P',D')\Big)\in {\cal I}\times {\cal I}: \ P\ne P',D=D' \Bigg\} \nn\\
{\cal O}_3 &=& \Bigg\{\Big((P,D),(P',D')\Big)\in {\cal I}\times {\cal I}:  \ P\ne P',D\ne D', P \in D'\Bigg\} \nn\\
{\cal O}'_3 &=& \Bigg\{\Big((P,D),(P',D')\Big)\in {\cal I}\times {\cal I}:  \ P\ne P',D\ne D', P' \in D\Bigg\} \nn\\
{\cal O}_4 &=& \Bigg\{\Big((P,D),(P',D')\Big)\in {\cal I}\times {\cal I}:  \ P\ne P',D\ne D',  P'\not\in D,P'\not\in D\Bigg\} \ .  \nn
\eeqa
These sets contain respectively $21, 42, 42, 84,84$ and $168$ elements.
\end{lemma}
\begin{demo}

$\bullet  \; {\cal D}$:  we shall show that the stabiliser of a point in ${\cal D}$ is of order 8.
Let $S \subset $Aut$(\F)$ be the stabiliser of $\Big(P,D\Big) \in \F$.
There is an  sequence of group homomorphisms 
\beqa
1 \to \text{Ker}(\Phi) \to S \to S_2\to 1 \nn
\eeqa
where $S_2$ is identified with the subgroup of permutations of $D$ which fix $P$ and $\Phi(s) \in S_2$ is the induced action
of $s \in S$ on $D$. By Proposition \ref{cor:4} (2),   $\text{Ker}(\Phi)$ is of order four and by Proposition \ref{cor:4} (4),
$\Phi$ is surjective (the complement of  $\text{Ker}(\Phi)$ in $S$ consists of two elements of order two and two elements of order four).
Hence $S$ is of order eight and $\A(\F)/S$ has $21=168/8$ elements. 
A straightforward count shows that ${\cal D}$ has  $21$ elements
 and it follows that ${\cal D}$ is a single orbit of  Aut$(\F)$.

$\bullet  \; {\cal O}_2$:   we shall show that the stabiliser $S$ of a point in ${\cal O}_2$ of order four.
Let $\Big((P,D),(P',D)\Big)\in {\cal O}_2$ and
suppose $g \in$Aut$(\F)$  fixes $\Big((P,D),(P',D)\Big)$, {\it i.e.},
\beqa
g\cdot P =P \ ,  \ \  g\cdot P' =P' \ , \ \ g\cdot D =D \  \nn \ .
\eeqa
Since $g$ fixes $P$ and $P'$ its fixes $P''$ the third point of the line $D$, and hence by  Proposition \ref{cor:4} since it fixes three points, $g$ is either of order one or two.
In fact by  Proposition \ref{cor:4} (2), either $g$ is the identity or $g$  can be taken to permute two pairs of points of the complement of $D$. There are only three elements of
this type and together with the identity it is clear that they form a subgroup of order four.
Hence $\A(\F)/S$ has $42=168/4$ elements. 
A straightforward count shows that ${\cal O}_2$ has  $42$ elements
 and it follows that ${\cal O}_2$ is a single orbit of  Aut$(\F)$.
Similarly ${\cal O}_1$ is a single orbit of $\A(\F)$ containing 42 elements.

$\bullet  \; {\cal O}_3$: we shall show that the stabiliser  $S$ of a point in ${\cal O}_3$ is $\mathbb S_2$. Let $\Big((P,D),(P',D')\Big)\in {\cal O}_3$ and
suppose $g \in$Aut$(\F)$  fixes $\Big((P,D),(P',D')\Big)$, {\it i.e.},
\beqa
g\cdot P =P \ ,  \ \ g\cdot D =D \ ,\ \  g\cdot P' =P' \ , \ \  g\cdot D' =D' \nn \ .
\eeqa
Since $g$ fixes $P$ and $P'$ its fixes $P''$ the third point of the line $D'=P\wedge P'$ and hence by Proposition \ref{cor:4}, since it fixes three points, $g$ is either of order one or two.
Then, either $g$ fixes each point of $D$ in which case  ({\it cf} Proposition \ref{cor:4}) $g$ is  the identity,  or $g$ permutes the two points of $D\setminus\{P\}$ and also
the two points $\F\setminus\big(D\cup D'\big)$.
Hence $\A(\F)/S$ has $84=168/2$ elements. 
A straightforward count shows that ${\cal O}_3$ has  $84$ elements
 and it follows that ${\cal O}_3$ is a single orbit of  Aut$(\F)$.
Similarly ${\cal O}'_3$ is a single orbit of $\A(\F)$ containing 84 elements.

$\bullet  \; {\cal O}_4$: we shall show that the stabiliser $S$ of an element of ${\cal O}_4$ is trivial. First let $\Big((P,D),(P',D')\Big)\in {\cal O}_4$ and
suppose $g \in$Aut$(\F)$  fixes $\Big((P,D),(P',D')\Big)$, {\it i.e.},
\beqa
g\cdot P =P \ ,  \ \ g\cdot D =D \ ,\ \  g\cdot P' =P' \ , \ \  g\cdot D' =D' \nn \ .
\eeqa
Since $g$ fixes $P$ and $P'$ its fixes $P''$ the third point of the line $P\wedge P'$.
Similarly, since  $g$  fixes $D$ and $D'$, it fixes the point $D\cap D'$ which is distinct from $P$ and $P'$  ($P \not \in D'$ and $P' \not \in D$).
The point $D\cap D'$ is also distinct from $P''$ since $D \neq D'$. We have now shown that $g$ fixes four distinct points of $\F$ and it  follows from Proposition \ref{cor:4}\
that $g$ is the identity. 
Hence $\A(\F)/S$ has $168$ elements. 
A straightforward count shows that ${\cal O}_4$ has  $168$ elements
 and it follows that ${\cal O}_4$ is a single orbit of  Aut$(\F)$.
\end{demo}

\begin{theorem}\label{theo:g2CR}
Let $\F$ be an Fano plane and let $\e$ be a composition factor.
Let $P, P'$ be disctinct point in $\F$  and let  $D, D'$ be  distinct lines in the Fano plane.
\begin{enumerate}
\item If $\Big((P,D),(P,D') \Big)\in {\cal O}_1$  then
\beqa
[X_{P,D},X_{P,D'}]=0 \ . \nn 
\eeqa
\item If $\Big((P,D),(P',D) \Big)\in {\cal O}_2$  then
\beqa
[X_{P,D},X_{P',D}]=2 \epsilon_{PP'} X_{P+P',P\wedge P'} \ . \nn 
\eeqa
\item If $\Big((P,D),(P',D') \Big)\in {\cal O}_3\cup {\cal O}'_3$ ({\it i.e.}, $P\in D',P'\not \in D$ or $P\not \in D', P'\in D$)  then
%$p' \in D$ or $p \in D'$ then
\beqa
[X_{P,D},X_{P',D'}]=- \epsilon_{PP'}  X_{P+P',P\wedge P'} \ . \nn 
\eeqa
\item  If $\Big((P,D),(P',D') \Big)\in {\cal O}_4$  ({\it i.e.}, $P\not\in D',P'\not \in D $)  then 
\beqa
[X_{P,D},X_{P',D'}]=- \epsilon_{PP'}  X_{P+ {P'},D + D'} \ ,\nn
\eeqa
where $D + D'$ is the third line containing  the point $D \wedge D'$.
\end{enumerate}
\end{theorem}

\begin{demo}
Parts (1)  of the theorem have already been proved   see Definition/Proposition \ref{DP} (3)) and part (2) is a straightforward calculation using \eqref{eq:so7}.
To prove the rest of  the theorem we first give the transformation law of brackets of two $X_{P,D}$
under the action of Aut$(\hat \F_\epsilon)$.    
We then explicitly calculate the bracket corresponding to  particular elements of each of the orbits ${\cal O}_i$ and then, using the transformation law
we deduce the formul\ae \ for  brackets in the general case. 

\begin{lemma}
Let $(P,D), (P',D') \in {\cal I} \times {\cal I}$  be such that $P \ne P'$. Suppose that
\beqa
\label{eq:brg2}
[X_{P,D},X_{P',D'}] = \alpha X_{P+ P',D''}  
\eeqa
where $\alpha \in \mathbb F$,  and $D''\in \F^\ast$. Then for all $g \in \mathrm{Aut}(\F)$
\beqa
 [X_{g\cdot P,g\cdot D},X_{g\cdot P',g\cdot D'}]&=&\alpha  \delta^\ast( g,P \wedge P') X_{g\cdot(P+ P'),g\cdot D''} \ , \nn
\eeqa
where $P\wedge P'$ is the line passing through $P$ and $P'$.
\end{lemma}
\begin{demo}
Let $\hat g \in \mbox{Aut}(\hat \F_\epsilon)$ be such that $\pi(\hat g) = g$ (see Proposition \ref{prop:tL} (2)). Recall that ({\it cf} \  Section \ref{sec:Auth}) we have
\beqa
\hat g X_{Q,L} \hat g^{-1} = \delta(\hat g,Q) X_{g\cdot Q,g\cdot L} \ \ \ \ \forall \; (Q,L) \in {\cal I} \ .\nn 
\eeqa
Hence conjugating  \eqref{eq:brg2} by $\hat g$  we get
\beqa
 \delta( \hat g,P) \delta( \hat g,P')[X_{g\cdot P,g\cdot D},X_{g\cdot P',g\cdot D'}]&=&\alpha\delta( \hat g,P+P')X_{g\cdot(P+P'),g\cdot D''} \ ,  \nn
\eeqa
which can be rewritten
\beqa
[X_{g\cdot P,g\cdot D},X_{g\cdot P',g\cdot D'}]&=&\alpha  \delta( \hat g,P) \delta( \hat g,P') \delta( \hat g,P+P')X_{g\cdot(P+P'),g\cdot D''} \ . \nn
\eeqa
However using Proposition \ref{prop:espdel} this reduces to
\beqa
[X_{g\cdot P,g\cdot D},X_{g\cdot P',g\cdot D'}] = \alpha \delta^\ast(g,(P,P')) X_{g\cdot (P+ P'),g\cdot D''} \ .\nn
\eeqa
\end{demo}

To choose a particular element of ${\cal O}_3$ we use the notation of  Figure \ref{fig:OF}.
Consider $\Big((P=P_1,D=D_1), (P'=P_3,D'=D_7)\Big) \in {\cal O}_3$ in which case
$(P+P''=P_7,D''=D_7)$.  Using \eqref{eq:so7} and \eqref{eq:X} we obtain by explicit calculation
\beqa
[X_{P_1,D_1},X_{P_3,D_7}]=-X_{P_7,D_7} \ , \nn
\eeqa
and by the lemma above this proves (3) since any element of ${\cal O}_3$ is conjugate to $\Big((P_1,D_1), (P_3,D_7)\Big)$. By symmetry this also proves (4).

To prove (5) consider  $\Big((P=P_4,D=D_1), (P'=P_5,D'=D_2)\Big) \in {\cal O}_4$ in which case
$(P+P'=P_7,D+D' =D_7)$.   Using \eqref{eq:so7} and \eqref{eq:X} we obtain  by explicit calculation
\beqa
[X_{P_4,D_1},X_{P_5,D_2}]=-X_{P_7,D_6} \ , \nn
\eeqa
and by the lemma above  this proves (5) since any element of ${\cal O}_4$ is conjugate to  $\Big((P_4,D_1), (P_5,D_2)\Big)$.
\end{demo}
\begin{corollary}\label{prop:g2}
%{\color{white}toto}
Let $(\F,\tau)$ be an oriented Fano plane and let  $\e$ be the canonical composition factor.
\begin{enumerate}
\item There is a decomposition as a direct sum of Cartan subalgebras:
\beqa
\mathfrak{g}_2(\mathbb F) = \bigoplus\limits_{P\in \F} \mathfrak{h}_{P} \ . \nn
\eeqa
\item If $P, Q, R$ are  non-aligned   then $\mathfrak{h}_{P} \oplus \mathfrak{h}_{Q} \oplus \mathfrak{h}_{R}$ generates $\g_2(\mathbb F)$.
\item If $P, Q, R$ are  three distinct points on the line $D$ such that  $\epsilon_{PQ} =\epsilon_{QR}= \epsilon_{RP}=1$ then:
\beqa
\big[X_{P,D}, X_{Q,D}\big] &=& \phantom{-} 2 X_{R,D} \ ,  \nn \\ 
\big[Y_{P,D}, Y_{Q,D}\big] &=&-2 Y_{R,D}  \ , \nn\\
\big[X_{S_1,D}, Y_{S_2,D}\big] &=& 0 \ , \ \  \forall S_1,S_2 \in D \nn \ .
\eeqa
In particular:
\begin{enumerate}
\item[$\bullet$]  $[\mathfrak{h}_P,\mathfrak{h}_Q]=\mathfrak{h}_R$.

\item[$\bullet$]
$
\mathfrak{g}_D=\mathfrak{h}_{P} \oplus \mathfrak{h}_{Q} \oplus \mathfrak{h}_{R}
$
is a Lie subalgebra of $\mathfrak{g}_2(\mathbb F)$ isomorphic to $\mathfrak{so}(4,\mathbb F)$.
\item[$\bullet$] If we define
\beqa
{\cal I}_{X,D}&=& \mathrm{Vect}\Big<X_{P,D}, X_{Q,D},X_{R,D}\Big> \ ,  \nn\\
{\cal I}_{Y,D}&=& \mathrm{Vect}\Big<Y_{P,D}, Y_{Q,D},Y_{R,D}\Big> \ ,  \nn\
\eeqa
then ${\cal I}_{X,D}$ and ${\cal I}_{Y,D}$ are the two ideals of $\mathfrak{g}_D$ each isomorphic to $\mathfrak{so}(3,\mathbb F)$.

\item[$\bullet$] The decomposition
\beqa
\mbox{Im} (\mathbb O_\F) = \mbox{Vect}\Big< e_P, \ P \in D \Big> \oplus  \mbox{Vect}\Big< e_P,\  P \not \in D \Big> \ , \nn
\eeqa
is stable under the action of  $\mathfrak{so}(4,\FF)$ and $ \mbox{Vect}\Big< e_P,\  P  \not \in D \Big> $ is isomorphic to the defining representation.
The representation $\mbox{Vect}\Big< e_P, \ P \in D \Big>$ is trivial for the action of ${\cal I}_{X,D}$ and isomorphic to the adjoint representation of ${\cal I}_{Y,D}$.
\end{enumerate}
\end{enumerate}
\end{corollary}

%\begin{demo}
%Straightforward except for the identification of the representations in part three. This can  best be  done by simply writing the matrices of the $X_{P,D}, Y_{Q,D}$ in the basis
%$\big<e_{P_i}, i =1,\cdots,7\big>$.
%\end{demo}

\begin{demo}
All parts of the Corollary are direct consequences of the Theorem above.
\end{demo}

\begin{remark}
As a consequence of this proposition one can associate to each line  of the oriented Fano plane an $\mathfrak{so}(4,\mathbb F)$ Lie subalgebra of $\g_2(\mathbb F)$ just as
by Proposition \ref{prop:oct} one can associate to each line of the oriented Fano plane a quaternion subalgebra of $\mathbb O_\F$.
\end{remark}

We now show that
one can associate a rank two simple Lie algebra to each point of the oriented Fano plane.

\begin{corollary}
For any point $P \in \F$, the Lie algebra $$ \mathfrak{s}_P={\rm{Vect}}\big\{X_{Q,D} \in \g_2(\FF) \ {\rm{s.t.}} \ X_{Q,D} \cdot e_P=0\big\}$$ is
isomorphic to $\mathfrak{sl}(3,\FF)$ if $\sqrt{-1} \in \FF$ and to $\mathfrak{su}(3,\mathbb F)$ if $\sqrt{-1} \not\in \FF$
(see  the beginning of Section  \ref{sec:cv-19} for the definition of $\mathfrak{su}(3,\mathbb F)$).
\end{corollary}

\begin{demo}
By Definition/Proposition \ref{DP} (1) if $P \in D$ then $[X_{Q,D},e_P]=0$.  Using the notation of Figure \ref{fig:OF}, without loss of generality we can suppose $P=P_1$.
Thus  $[X_{Q,D},e_{P_1}]=0$ if
\beqa
X_{Q,D}= \left\{
\begin{array}{lll}
 X_{P_1,D_1}\ ,&X_{P_2,D_1}\ ,&X_{P_4,D_1}\ , \\
 X_{P_7,D_7}\ ,&X_{P_1,D_7}\ ,&X_{P_3,D_7}\ , \\
 X_{P_5,D_5}\ ,&X_{P_6,D_5}\ ,&X_{P_1,D_5}\ ,
\end{array}
\right. \nn
\eeqa
with $X_{P_1,D_1}+X_{P_1,D_7}+X_{P_1,D_5}=0$ (see Definition/Proposition \ref{DP} (3)). Thus we have eight linearly independent generators
which clearly span a Lie subalgebra of $\mathfrak{g}_2(\FF)$.
Again from  Definition/Proposition \ref{DP} (4), $X_{P_1,D_1}, X_{P_1,D_7}$ span the Cartan subalgebra $\mathfrak h_{P_1}$.
We now introduce (in $\mathfrak g_2(\mathbb F)$ if $\sqrt{-1} \in \FF$ or in $\mathfrak g_2\Big(\mathbb F[\sqrt{-1}]\Big)$ if $\sqrt{-1}\not  \in \FF$) the elements:
\beqa
h_{P_1,D_1} = -\sqrt{-1} X_{P_1,D_1} \ , \nn\\
h_{P_1,D_7} = -\sqrt{-1} X_{P_1,D_7} \ ,\nn
\eeqa
which generate  a Cartan subalgebra, together with the elements
\beqa
e^\pm_{P_1,D_1}&=& X_{P_2,D_1}\mp \sqrt{-1} X_{P_4,D_1} \ , \nn\\
e^\pm_{P_1,D_7}&=& X_{P_3,D_7}\mp \sqrt{-1} X_{P_7,D_7} \ . \nn
\eeqa
Then from Theorem \ref{theo:g2CR} we have
\beqa
\begin{array}{lllllll}
\Big[h_1,e^\pm_{P_1,D_1}\Big]&=&\pm 2 e^\pm_{P_1,D_1}\ , &&\Big[h_1,e^\pm_{P_1,D_7}\Big]&=&\mp  e^\pm_{P_1,D_7}\ , \\
\Big[h_2,e^\pm_{P_1,D_1}\Big]&=&\mp  e^\pm_{P_1,D_1}\ , &&\Big[h_2,e^\pm_{P_1,D_7}\Big]&=&\pm 2  e^\pm_{P_1,D_7}\ , 
\end{array}\nn
\eeqa
and
\beqa
\Big[e^+_{P_1,D},e^-_{P_1,D'}\Big]&=& \left\{
\begin{array}{cll}
-4 h_{P_1,D}& \mathrm{if}&D=D' \ , \nn\\
0&\mathrm{if}&D\ne D' \ .
\end{array}
\right.
\eeqa
Thus the Cartan matrix is given by
\beqa
A =\begin{pmatrix} \phantom{-}2&-1\\-1&\phantom{-}2\end{pmatrix} \ . \nn
\eeqa
The last two generators of $\mathfrak{s}_P$ are given by
\beqa
[e^\pm_{P_1,D_1},e^\pm_{P_1,D_7}]=-2 e^\pm_{P_1,D_5}= -2(X_{P_5,D_5} \pm \sqrt{-1} X_{P_7,D_5}) \ , \nn
\eeqa
and satisfy (again by Theorem \ref{theo:g2CR})
\beqa
\begin{array}{lllllll}
\Big[h_1,e^\pm_{P_1,D_5}\Big]&=&\pm e^\pm_{1,D_5} \ , &&\Big[h_2,e^\pm_{P_1,D_5}\Big]&=&\pm  e^\pm_{P_1,D_5}\ ,
\end{array}\nn
\eeqa
and
\beqa
\Big[e^+_{P_1,D_5},e^-_{P_1,D_5}\Big]&=& -4(h_{P_1,D_1}+h_{P_1,D_7})\ .\nn
\eeqa
This ends the proof.
\end{demo}

\begin{remark}
The Lie algebra $\mathfrak{s}_P$ has a natural representation on the vector space:
\beqa
V= {\rm{Vect}}\big\{e_Q \ {\rm s.t.} \  Q \in \F\ , Q \ne P \big\}\ . \nn
\eeqa
This six-dimensional vector space is stable under ${\mathfrak s} _P$ and carries a natural
 $\mathfrak{s}_P-$invariant ``almost-complex structure'' $J \in {\rm{End}}(V)$ defined by
\beqa
J(e_Q)=\epsilon_{QP} \; e_{P+Q} \ . 
                  \nn
\eeqa
It also carries the ${\mathfrak s}_P-$invariant non-degenerate quadratic form
\beqa
g(e_Q,e_R) = \delta_{QR} \nn
\eeqa
for which, one checks, $J$ is an isometry. The Lie algebra $\mathfrak{s}_P$ is then
characterised as  the Lie subalgebra of ${\rm End}(V)$ which
preserves $g$ and $J$. This  representation is irreducible iff $\sqrt{-1} \not \in \FF$.
\end{remark}

%%%%%%%%%%%%%%%%%%%%%

%%%%%%%%%%%%%%%%%%%%%%%%%%%%%%%%%%%%%%%
Proposition  \ref{prop:oct} associates  composition subalgebras of the octonions to certain configurations of points in the Fano plane.
Analogously, one can  associate  Lie subalgebras of $\g_2(\FF)$ to elements of ${\cal I}\times {\cal I}$ in the following way:

\begin{itemize}
\item[${\cal O}_1$:] For each element $\Big((P,D),(P,D')\Big)$ of ${\cal O}_1$  the  Lie subalgebra of $\mathfrak{g}_2(\mathbb F)$ generated by $X_{P,D}, X_{P,D'}$
is a Cartan subalgebra denoted $\mathfrak{h}_P$ in Proposition \ref{DP}.
\item[${\cal O}_2$:] For  each element $\Big((P,D),(P',D)\Big)$ of ${\cal O}_2$  the  Lie subalgebra of $\mathfrak{g}_2(\mathbb F)$
generated by $X_{P,D},X_{P',D}$  is isomorphic to $\mathfrak{so}(3,\mathbb F)$ by  Proposition \eqref{prop:g2} (2).
\item[${\cal O}_3$:] For each element  $\Big((P,D),(P',D')\Big)$ of ${\cal O}_3$  the Lie subalgebra of $\mathfrak{g}_2(\mathbb F)$, generated by $X_{P,D},X_{P',D'}$ is isomorphic to
 $\mathfrak{so}(3,\mathbb F) \times \mathfrak{so}(2,\mathbb F))$. For instance, in the case of $X_{P_1,D_1}$ and $X_{P_7,D_7}$  we have: 
 %$X_{p,D},X_{p',D'}$ satisfies: % and this turns out  isomorphic to  $\mathfrak{so}(3,\mathbb F) \rtimes \mathfrak{so}(2,\mathbb F))$. For instance
\beqa
\big[X_{P_1,D_1},X_{P_7,D_7}\big]&=& X_{P_3,D_7} \nn\\
\big[X_{P_1,D_1},X_{P_3,D_7}\big]&=&- X_{P_7,D_7} \nn\\
\big[ X_{P_1,D_1}, X_{P_1,D_7} \big]&=&0\nn \ . 
\eeqa
It follows that $X_{P_1,D_1} + \frac12 X_{P_1,D_7}= -\frac12 Y_{P_1,D_7}$ is in the centre of this Lie algebra which  for brevity we  shall denote $\mathfrak{s}$. The
$X_{P_7,D_7}, X_{P_3,D_7}, X_{P_1,D_7}$ generate an ideal  of $\mathfrak{s}$  isomorphic to $ \mathfrak{so}(3,\mathbb F)$. 
The Lie subalgebra $\mathfrak{s}$ can also be  obtained by adding an element of one of the ideals of the $ \mathfrak{so}(4,\mathbb F)$ associated to $D_7$ to  the other ideal (see Proposition \ref{prop:g2} (2)).
Under the action of  $\mathfrak{s}$  we have the decompositions:
\beqa
\mbox{Im} (\mathbb O_\F)  &=&\mbox{Vect}\Big< e_P, \ P\in D_7 \Big> \oplus  \mbox{Vect}\Big< e_P,\  P  \not \in D_7 \Big> \nn\\
&=&\Big<e_{P_7},e_{P_1},e_{P_3}\Big> \oplus\Big<e_{P_2},e_{P_4},e_{P_5},e_{P_6}\Big> \nn\\
&=&
\left\{
\begin{array}{ll}
 \Big({\bf 1_0} \oplus {\bf 2}\Big) \oplus {\bf 4} & \text{if~} -1 \ \text{is~not~a~square~in~} \FF\\
 \Big({\bf 1_0} \oplus {\bf 1} \oplus {\bf 1'}\Big) \oplus\Big( {\bf 2 }\oplus {\bf 2'}\Big)  & \text{if~} -1 \ \text{is~a~square~in~} \FF\\
\end{array}
\right.\nn
\eeqa
where ${\bf 1_0}$ is the trivial one-dimensional representation, ${\bf 1}, {\bf 1'}$ are irreducible one-dimensional representations,  ${\bf 2}, {\bf 2'}$ are irreducible two-dimensional representations
and ${\bf 4}$ is an irreducible four-dimensional representation.
\item[${\cal O}'_3$:] Same as ${\cal O}_3$.
\item[${\cal O}_4$:] Finally for  each element  $\Big((P,D),(P',D')\Big)$ of  ${\cal O}_4$ the Lie algebra generated by $X_{P,D}, X_{P',D'}$ is  $\g_2(\mathbb F)$.
\end{itemize}

%\phantomsection
%\addcontentsline{toc}{chapter}{Bibliographie}
\bibliographystyle{utphys}
\bibliography{ref}

\end{document}